\let\ORIlabel\label
\let\ORIrefstepcounter\refstepcounter
\AddToHook{package/hyperref/before}{
   \let\label\ORIlabel 
   \let\refstepcounter\ORIrefstepcounter}
\documentclass[hidelinks,onefignum,onetabnum]{siamart220329}

\usepackage{lipsum}
\usepackage{amsfonts}
\usepackage{graphicx}
\usepackage{epstopdf}
\usepackage{algorithmic}
\usepackage{multicol}
\usepackage{multirow}
\usepackage{accents}
\ifpdf
  \DeclareGraphicsExtensions{.eps,.pdf,.png,.jpg}
\else
  \DeclareGraphicsExtensions{.eps}
\fi

\usepackage{multirow} 
\usepackage{makecell} 
\usepackage{tikz}
\usepackage{pgfplots}
\usetikzlibrary{positioning}
\usetikzlibrary{plotmarks}
\usetikzlibrary{patterns}
\newlength\figureheight
\newlength\figurewidth

\usepackage{booktabs}
\usepackage{amssymb}
\usepackage{xcolor}
\usepackage[table]{xcolor}
\usepackage[dvipsnames]{xcolor}
\newcommand{\R}{\mathbb{R}}

\newcommand{\matrx}[1]{\mathbf{#1}}
\renewcommand{\vec}[1]{\mathbf{#1}}
\newcommand{\A}{\matrx{A}} 

\newcommand{\store}{\mathrm{R}}
\newcommand{\solve}{\mathrm{S}}

\newcommand{\TG}{\mathrm{TG}}
\newcommand{\V}{\mathrm{V}}
\newcommand{\C}{\mathrm{C}}
\newcommand{\IC}{\mathrm{IC}}

\definecolor{Green}{RGB}{0, 128, 0}
\definecolor{_LimeGreen}{RGB}{232, 242, 161}
\newcommand{\tpc}[1]{\textcolor{black}{#1}}
\newcommand{\vpc}[1]{\textcolor{black}{#1}}
\newcommand{\tp}{\tpc{\dot{\varepsilon}}}
\newcommand{\vp}{\vpc{\dot{\varepsilon_j}}}

\newsiamremark{remark}{Remark}
\newsiamremark{hypothesis}{Hypothesis}
\crefname{hypothesis}{Hypothesis}{Hypotheses}
\newsiamthm{claim}{Claim}

\headers{Mixed precision multigrid with IC smoothing}{P. Vacek, H. Anzt, E. Carson, N. Kohl, U. Rüde, Y.-H. Tsai}

\title{Mixed precision multigrid with smoothing based on incomplete Cholesky factorization\thanks{Submitted to the editors November 7th, 2025. Revised version submitted on July 26th 2026.
\funding{The first and third authors were supported by the European Union (ERC, inEXASCALE, 101075632). The first and fifth author were also supported by the EuroHPC JU grant agreement 101144014. 
The first author would also like to thank for the support from a national grant managed by the French National Research Agency (Agence Nationale de la Recherche) attributed to the Exa-SofT project of the NumPEx PEPR program part of the ”France 2030” initiative, under the reference ANR-22-EXNU-0003.
The third author additionally acknowledges support from the Charles University Research Centre program, UNCE/24/SCI/005.
The fifth author received financial support of the European Union under the REFRESH -- Research Excellence For Region Sustainability
and High-tech Industries project number CZ.10.03.01/00/22\_003/0000048
via the Operational Programme Just Transition.
Views and opinions expressed are those of the authors only and do not necessarily reflect those of the European Union, the European Research Council or the EuroHPC JU. Neither the European Union nor the granting authority can be held responsible for them.
}}}

\author{Petr Vacek\thanks{CNRS-IRIT, Toulouse (\email{petr.vacek@irit.fr}), IFP Energies Nouvelles, Rueil-Malmaison, France, and Department of Numerical Mathematics, Charles University, Prague, Czech Republic.} 
\and Hartwig Anzt\thanks{Technical University of Munich - Campus Heilbronn, Heilbronn,
Germany (\email{\{hartwig.anzt,yu-hsiang.tsai\}@tum.de}).}
\and Erin Carson \thanks{Department of Numerical Mathematics, Charles University, Prague, Czech Republic  (\email{carson@karlin.mff.cuni.cz}).}
\and Nils Kohl\thanks{Ludwig-Maximilians-Universität München - Department of Earth and Environmental Sciences, Munich, Germany (\email{nils.kohl@lmu.de}).}
\and Ulrich Rüde\thanks{Department of Computer Science, FAU Erlangen-Nürnberg, Germany; Department of Applied Mathematics, VSB-Technical University of Ostrava,
Ostrava, Czech Republic;  Centre Européen de Recherche et de Formation Avancée en Calcul Scientifique,
France (\email{ulrich.ruede@fau.de}).}
\and Yu-Hsiang Tsai \footnotemark[3]}
\usepackage{amsopn}

\ifpdf
\hypersetup{
  pdftitle={Mixed precision multigrid with smoothing based on incomplete Cholesky factorization},
  pdfauthor={P. Vacek, Y-H. Tsai, H. Anzt, U. Rüde, E. Carson}
}
\fi

\begin{document}
\maketitle
\begin{abstract}
Multigrid methods are popular iterative methods for solving large-scale sparse systems of linear equations.
We present a mixed precision formulation of the multigrid V-cycle method with general assumptions on the finite precision errors coming from the application of coarsest-level solver and smoothing.
Inspired by existing analysis, we derive a bound on the relative finite precision error of the V-cycle which gives insight into how the finite precision errors from the individual components of the method may affect the overall finite precision error.
We use the result to study V-cycle methods with smoothing based on incomplete Cholesky factorization.
The results imply that in certain settings, the precisions used for applying the incomplete Cholesky smoothing can be significantly lower than the precision used for computing the residual, restriction, prolongation, and correction on a given level.
We perform numerical experiments using simulated floating point arithmetic with the MATLAB Advanpix toolbox as well as experiments on GPUs using the Ginkgo library.
The experiments illustrate the theoretical findings and show that in the considered settings, the incomplete Cholesky smoothing can be applied in relatively low precisions, resulting in significant improvements in the execution time (up to $28\%$ less),  energy savings (up to $26\%$ less) and required memory (up to $16\%$ less) in comparison with the uniform double precision variant.
\end{abstract}



\section{Introduction}
There is extensive ongoing research in numerical methods capable of exploiting multiple precisions; see, e.g., the surveys \cite{Higham2022,Abdelfattah2021}.
In some cases, such mixed precision methods can achieve the same overall accuracy as their uniform precision counterparts in a shorter amount of time, requiring less memory, and consuming less energy.
In this text we study mixed precision variants of multigrid methods \cite{Trottenberg2001,Briggs2000,Brandt2011} which are frequently used when solving systems of linear equations. 
Multigrid methods can be applied both as standalone solvers and as preconditioners for iterative methods. 
The computation relies on having a hierarchy of problems, which can be obtained either by discretizing a continuous problem on multiple nested meshes (geometric multigrid) or constructed based on the properties of the system matrix (algebraic multigrid).
The approximate solution is computed using smoothing on fine levels and by solving a system of linear equations on the coarsest level. 
Smoothing on any level should contribute primarily to reducing the high frequency components of the error, while the low frequency components are eliminated via a coarse grid correction.
There are different multigrid schemes (V-cycle, W-cycle, full multigrid) that vary in the pattern in which the individual levels are visited during the computation.

Implementations of multigrid methods which employ different precision formats in different parts of the method
have been developed and tested on various problems; see, e.g., \cite{Tsai2023,Tsai2023b,Tsai2024,Zong2024}.
The first finite precision error analysis of mixed precision multigrid methods was presented in \cite{McCormick2021} and further extended in \cite{McCormick2023}. The results were used by the authors in a paper focusing on achieving discretization error accuracy when solving elliptic PDEs \cite{Tamstorf2021} and adapted also for multigrid methods with block floating point arithmetic in \cite{Kohl2023}.

The finite precision error analysis of the V-cycle method presented in \cite{McCormick2021,McCormick2023} is based on viewing the method as an iterative refinement (IR) process on the finest level with a V-cycle starting with zero initial approximation as the inner solver.
This point of view enables separation of the analysis into the analysis of IR and the analysis of one V-cycle with zero initial approximation. 
The authors consider the case where the computations on different levels in the V-cycle are potentially performed in different precisions with different unit roundoffs. The operations on a concrete fine level, i.e., computation of the residual, restriction, prolongation, correction, and smoothing, are assumed to be all done in the same precision. 
The results imply that in certain settings the precisions used on the coarse levels can be chosen progressively lower and lower without having a significant effect on the convergence rate.

Multigrid methods are, in practice, also applied with computationally intensive smoothers.
Smoothing routines based on incomplete Cholesky (IC) or incomplete LU factorization are, for example, used when solving elliptic PDEs with large anisotropy and/or when using discretization based on high-degree polynomial basis functions; see, e.g., the early papers \cite{Kettler1982,Wesseling1982,Wesseling1982b,Kettler1986} or \cite{Tielen2020,Drzisga2023,Thomas2024}.
To use IC smoothing, the IC factorization must be precomputed. Each application then requires solving triangular systems with the IC factor and its transpose. 

In this work we ask whether the mixed precision approach could be also used to speed up the application of the smoothing routines, whose role is to reduce the high frequency components of the error. 
For IC smoothers, this could mean storing the IC factors in low-precision and/or solving the triangular systems in low-precision. 
This opens a series of questions. What precisions should be used in the mentioned stages of the IC smoother? How should these precisions be chosen with respect to the application of the smoother inside the V-cycle method? 

Motivated by these questions, we present a formulation and finite precision error analysis of the V-cycle method with general assumptions on the smoothers and the coarsest-level solver. Rather than assuming that the smoothers and the coarsest-level solver are applied in a certain precision, we impose assumptions on the finite precision errors resulting from their applications. This enables us to consider also mixed precision smoothers and coarsest-level solvers.
The derived bound on the finite precision error gives insight into how the finite precision errors from the individual parts of the V-cycle may affect the overall finite precision error. 
The results are valid for both geometric and algebraic multigrid provided that the matrices on the coarse levels satisfy the Galerkin condition.
We further formulate a mixed precision IC smoothing routine and present a bound on the finite precision error resulting from its application. We assume that the triangular problems are solved using substitution. 
The analysis does not include the finite precision errors occurring when computing the IC factorization.
The results imply that in certain settings the precision used for applying the IC smoothing could be significantly lower than the precision used for computing the residual, restriction, prolongation, and correction on a given level.
We test the theoretical results and performance of the presented methods through a series of numerical experiments. We solve systems coming from finite element (FE) discretization of the Poisson equation. We run experiments with simulated floating point arithmetics in MATLAB using the Advanpix toolbox \cite{advanpix} as well as experiments performed on GPUs using the Ginkgo library \cite{ginkgo,Ginkgo2024}.

The paper is organized as follows. In \Cref{sec:notation_fp}, we establish notation, present the standard rounding model, and state bounds on the finite precision errors in basic vector and matrix operations. 
A mixed precision two-grid cycle is presented in 
\Cref{sec:tg} together with its finite precision error analysis.
These results are generalized to a multigrid V-cycle in \Cref{sec:V-cycle}.
In \Cref{sec:IC_smoothing}, we present a mixed precision smoothing routine based on IC factorization and derive a bound on the finite precision errors occurring in its application. The results on the effects of finite precision errors in a V-cycle with IC smoothing are summarized in \Cref{sec:IC-V-cycle}. Basic scaling strategies that may help prevent overflow and underflow errors are discussed in \Cref{sec:scaling}.
Numerical experiments illustrating the results and performance gains are presented in \Cref{sec:num-exp-simul,sec:num-exp-gingko}. 

\section{Notation, finite precision arithmetic, and standard rounding model}\label{sec:notation_fp}
We consider all vectors and matrices in this paper to be real. We denote the Euclidean inner product as $\langle \cdot, \cdot \rangle$, and the Euclidean vector or matrix norm as $\| \cdot \|$.
For a symmetric positive definite (SPD) matrix $\A$, we denote the $\A$-norm of a vector $\vec{v}$ as  $\| \vec{v}\|_\A= \sqrt{\langle \A \vec{v}, \vec{v} \rangle}$; we use the same notation for the associated matrix norm. Below we use the following relations between the Euclidean and the $\A$ vector norms. For any vector $\vec{v}$ it holds that (see \Cref{ape:norms})
\begin{align}
\label{eq:A_norm_by_Euclid}
\| \vec{v} \|_\A & \leq \| \A \|^{\frac{1}{2}} \| \vec{v} \|, \\ 
\label{eq:Euclid_by_A_norm}
 \| \vec{v} \| &\leq \| \A^{-1} \|^{\frac{1}{2}} \| \vec{v} \|_{\A}, \\
\label{eq:A*v_by_A_norm}
\| \A \vec{v} \| & \leq \| \A \|^{\frac{1}{2}} \| \vec{v} \|_{\A}, \\
\label{eq:A^-1_v_by_Euclid}
 \| \A^{-1} \vec{v} \|_{\A} &\leq \| \A^{-1}\|^{\frac{1}{2}} \| \vec{v}\|.
\end{align}

For a matrix $\matrx{K}$ we denote by $| \matrx{K}|$ the matrix with the component-wise absolute values of the entries of $\matrx{K}$.
The condition number of an invertible matrix $\matrx{K}$ is denoted by $\kappa_{K}=\| \matrx{K}^{-1} \|  \| \matrx{K}\|$. By $\underline{\kappa}_{K}$ we denote a variant of the condition number containing $\||\matrx{K} | \|$ instead of $\|\matrx{K} \|$, i.e., $\underline{\kappa}_{K}=\| \matrx{K}^{-1} \| \| | \matrx{K} | \|$.
We use $\matrx{I}$ with various subscripts to denote identity matrices.

We consider the standard model for accounting for finite precision errors (see, e.g., \cite[Chapter 1]{Higham2002}), which is also used in the existing finite precision analysis of multigrid methods in \cite{McCormick2021,McCormick2023}. 
Consider a floating point arithmetic with unit roundoff $\varepsilon$. 
Rounding a vector $\vec{v}$ and a matrix $\matrx{K}$ to $\varepsilon$-precision results in, respectively, 
\begin{equation}\label{eq:fp_round}
    \vec{v} + \delta, \quad \| \delta \| \leq  \varepsilon \| \vec{v} \|, \quad \text{and} \quad  \vec{K} + \Delta \matrx{K}, \quad | \Delta \matrx{K} | \leq  \varepsilon | \matrx{K} |,
\end{equation}
where the second inequality is understood entry-wise.
Let $m_K$ denote the maximum number of nonzero entries in a row of matrix $\matrx{K}$ and let $\underline{m}_K$ denote the maximum number of nonzero entries in a row or a column of~$\matrx{K}$. 
Assume that $(m_K+2) \varepsilon < 1$ and that $\vec{v}$, $\vec{w}$ are vectors belonging to the $\varepsilon$-precision arithmetic. 
Computing $\vec{v} + \vec{w} $,  $\matrx{K} \vec{w}$  and $\vec{v} - \matrx{K}\vec{w}$, in $\varepsilon$-precision, where the matrix $\matrx{K}$ is first rounded to $\varepsilon$-precision, results in,  respectively,
\begin{align}
\label{eq:fp_add}
    &\vec{v} + \vec{w} + \delta, \quad \|\delta \| \leq  \varepsilon \| \vec{v} + \vec{w}\| , \\
\label{eq:fp_mat_vec_quant}
&\matrx{K}\vec{w} + \delta, \quad \|\delta \| \leq  
\frac{ (m_K+1) \varepsilon}{1 - (m_K+1) \varepsilon}
 \||\matrx{K} | \| \| \vec{w}\|,
\\
\label{eq:fp_mat_vec_add_quant}
&\vec{v} - \matrx{K}\vec{w} + \delta, \quad  \| \delta \|\leq 
\frac{ (m_K+2) \varepsilon}{1 - (m_K+2) \varepsilon}
 ( \|  \vec{v} \| +  \|| \matrx{K} | \| \| \vec{w}\|);
\end{align}
see e.g., \cite[Sections~2.2, 3.1 and 3.5]{Higham2002}.
Throughout the whole text, we assume that the computations do not break down due to overflow or underflow errors. This is a standard assumption in the literature.

\section{Two-grid method} 
\label{sec:tg}
In this section, we study the effects of finite precision errors in a two-grid method (TG) for solving $ \A\vec{y}=\vec{f}$, where $\A \in \R^{n \times n}$ is an SPD matrix and $\vec{f}\in \R^{n}$. 
We present a bound on the finite precision error after one TG cycle starting with zero initial approximation. The finite precision error analysis of multiple TG cycles can be then obtained by viewing the method as an iterative refinement method, with the TG cycle with zero initial approximation as the inner solver and using the results presented in \cite{McCormick2021, Tamstorf2021}.

In the TG cycle the approximate solution is computed using \emph{smoothing} and a \emph{coarse-grid correction}. For simplicity of the analysis we consider a version where smoothing is applied only before the coarse-grid correction; for other variants, see, e.g., \cite{Briggs2000,Trottenberg2001}. We first describe assumptions on the components of the method in exact arithmetic and then add the assumption on their application in finite precision.

We assume that for any approximation $\vec{v}$, the application of smoothing in exact arithmetic can be written $ \vec{v} + \matrx{M}(\vec{f}-\A \vec{v})$, where $\matrx{M} \in \R^{n \times n}$ is non-singular and 
\begin{equation}\label{eq:smoother_assumption_Anorm}
\| \matrx{I} - \matrx{M} \A \|_{\A}<1.
\end{equation}
Since we consider zero initial approximation, the above expression simplifies to $ \matrx{M}\vec{f}$.
The coarse-grid correction consists of computing the residual, restricting it to the coarse grid, solving the coarse-grid error equation, prolongating the correction to the fine grid, and correcting the previous approximation.
We assume that there exists an SPD coarse-grid matrix $\A_\C \in \R^{n_\C \times n_\C}$ and a full rank prolongation matrix $\matrx{P}\in \R^{n\times n_\C}$
such that the Galerkin condition is satisfied, i.e., $\A_{\C} = \matrx{P}^{\top}\A \matrx{P}$. 
We assume that the restriction matrix is the transpose of the prolongation matrix.
We further assume that for any vector $\vec{f}_{\C}$, the application of the coarse-grid solver in exact arithmetic can be written $\matrx{M}_\C \vec{f}_{\C}$,
where $\matrx{M}_\C$ is a non-singular matrix and 
\begin{equation}\label{eq:perturbed_coarse_solve_A_norm}
    \| \matrx{I}_\C-\matrx{M}_\C \A_{\C}\|_{\A_\C}<1.
\end{equation}
We note that for the exact coarse-grid solve, i.e., $\matrx{M}_\C = \A^{-1}_\C$,  \eqref{eq:perturbed_coarse_solve_A_norm} is automatically satisfied. The formulation here allows approximate linear coarse-grid solvers. 
 
We assume that the residual computation, prolongation, restriction, and correction are performed in arithmetic with a unit roundoff $\dot{\varepsilon}$, that the matrices $\A$ and $\matrx{P}$ are rounded to the $\dot{\varepsilon}$-precision before the computations, and that the unit roundoff $\dot{\varepsilon}$ is small enough such that $(m_A+2) \dot{\varepsilon} < 1$ and $(\underline{m}_P+1) \dot{\varepsilon} < 1$.

Rather than assuming that the smoother and coarse-grid solver are applied in a precision with a certain unit roundoff, we impose assumptions on the resulting relative finite precision errors. This allows us to consider also mixed precision smoothers or coarse-grid solvers.
We assume that there exist positive constants $\Lambda_M$ and $\Lambda_\C$
such that the application of the smoother and coarse-grid solver in finite precision for any vectors $\vec{f}$  and $\vec{f}_{\C}$ results in, respectively,
\begin{align}
\label{eq:smoother_assumption_fp}
    &\matrx{M}\vec{f} + \delta_M, \quad \| \delta_{M} \| \leq \Lambda_{M} \| \vec{f}\|, \\
    \label{eq:solver_assumption_fp}
 & \matrx{M}_\C \vec{f}_{\C} + \delta_\C, \quad \| \delta_{\C} \|_{\A_\C} \leq \Lambda_{\C} \| \A^{-1}_{\C} \vec{f}_{\C}\|_{\A_{\C}},    
\end{align}
and the resulting vectors belong to the $\dot{\varepsilon}$-precision arithmetic.

The TG cycle with zero initial approximation is formulated in \Cref{alg:two-grid}.
We assume that its application in exact arithmetic reduces the $\A$-norm of the error by a factor $\rho_{\TG}<1$, i.e.,   
\begin{equation}\label{eq:TG_exact_rho}
    \| \vec{y}- \vec{y}_{\TG} \|_{\A} \leq \rho_{\TG} \| \vec{y} \|_{\A},
\end{equation}
where $\vec{y}_{\TG}$ is the approximation computed by the TG cycle in exact arithmetic. For convergence analysis of multigrid methods in exact arithmetic see, e.g., \cite{Yserentant1993,Trottenberg2001,Notay2007,Notay2015}.

\begin{algorithm}[H]
\caption{Two-grid cycle with zero initial approximation $\mathbf{TG}(\vec{f})$.}
\label{alg:two-grid}
\begin{algorithmic}[1]
\STATE{$\vec{f} \leftarrow \mathrm{Round(\vec{f},\dot{\varepsilon}\text{-precision})}$} 
\STATE{$\vec{v}^{[1]} \leftarrow \matrx{M}\vec{f}$ \COMMENT{Apply smoothing.}}
\STATE{$\vec{r}^{[1]} \leftarrow \vec{f} -  
\A \vec{v}^{[1]}$ \COMMENT{Compute residual in $\tp$-precision.}}
\STATE{$\vec{r}^{[2]}_{\C} \leftarrow \matrx{P}^{\top} \vec{r}^{[1]}$ \COMMENT{Restrict the residual to the coarse grid in $\tp$-precision.}}
\STATE{$\vec{v}^{[2]}_{\C} \leftarrow \matrx{M}_\C \vec{r}^{[2]}_{\C}$
\COMMENT{Coarse-grid solve.
}}
\STATE{$\vec{v}^{[3]} \leftarrow \matrx{P} \vec{v}^{[2]}_{\C}$ \COMMENT{Prolongate the correction to the fine level in $\tp$-precision.}}
\STATE{$\vec{v}^{[4]} \leftarrow \vec{v}^{[1]} + \vec{v}^{[3]}$ \COMMENT{Correct the previous approximation in $\tp$-precision.}}
\RETURN{$\vec{v}^{[4]}$}
\end{algorithmic}
\end{algorithm}

We present the following result on the effects of finite precision errors in the TG cycle. Its proof can be found below.

\begin{theorem}\label{thm:tg} 
Let $\vec{y}_{\TG}$ and  $\hat{\vec{y}}_{\TG}$ be the approximate solution of $\A\vec{y}=\vec{f}$ computed using one TG cycle (\Cref{alg:two-grid}) applied exactly and in finite precision, respectively. 
The $\A$-norm of the finite precision error $\delta_{\TG}=\vec{y}_{\TG} - \hat{\vec{y}}_{\TG}$ can be bounded as 
\begin{equation}
\label{eq:thm:tg_difference}
\| \delta_{\TG}\|_{\A} \leq  \Lambda_{\TG} \| \vec{y} \|_{\A},
\quad
\Lambda_{\TG} =  \Lambda_\C + 3\| \A \| \Lambda_M 
  +  \dot{\varepsilon} \kappa^{\frac{1}{2}}_{A} (C_1  \|\matrx{M} \|
+ C_2) + R,
\end{equation}
where $C_1$ and $C_2$ are positive constants depending on   $\| \matrx{A}\|$, $ \| | \matrx{A} | \|$, $\| \matrx{P}\|$, $ \| | \matrx{P} | \|$, 
\begin{equation*}
    m_{A,\dot{\varepsilon}} =\frac{(m_A+2)}{1-(m_A+2)\dot{\varepsilon}}, \quad 
    \underline{m}_{P,\dot{\varepsilon}} =\frac{(\underline{m}_P+1)}{1-(\underline{m}_P+1)\dot{\varepsilon}}, \quad \text{and} \quad
\xi =  \frac{\| \A^{-1}_{\C}\|^{\frac{1}{2}}}{ \| \A^{-1} \|^{\frac{1}{2}}}.
\end{equation*}
The remainder $R$ contains additional higher order terms. 
If the sum of the first three terms in $\Lambda_{\TG}$ is sufficiently smaller than one, the remainder $R$ is negligible. 
If $\rho_{\TG} + \Lambda_{\TG} < 1$, the TG cycle applied in finite precision reduces the relative $\A$-norm of the error by at least the factor $ \rho_{\TG} + \Lambda_{\TG} $, i.e.,
$\| \vec{y} - \hat{\vec{y}}_{\TG}  \|_{\A}  \leq (\rho_{\TG} +  \Lambda_{\TG}) \| \vec{y} \|_{\A}$.
In other words, the TG cycle applied in finite precision is a contraction with a contraction factor $ \rho_{\TG} + \Lambda_{\TG}$.
\end{theorem}

We intentionally present this theorem without any additional assumptions on the individual terms in the estimates so that it can be used in various settings. A more detailed expression can be found inside the proof below if needed.

An important feature of the bound is that it provides insight into how the finite precision errors coming from applying the smoother or the coarse-grid solver may affect the overall finite precision error.
In particular, we see that the bound on the relative finite precision error of the coarsest-grid solver $\Lambda_\C$ is present as a standalone term.
The bound on the relative finite precision error of the smoother is multiplied by $3 \| \A \|$, but not for example by $\| \A^{-1} \|^{\frac{1}{2}}$ or $\kappa^{\frac{1}{2}}_{A}$, which might be large.
Another useful observation is that the choice of the smoother may affect the requirements on the $\dot{\varepsilon}$-precision through the term $\| \matrx{M}\|$.

We generalize the results to a multigrid V-cycle scheme in the next section.

\begin{proof}[Proof of \Cref{thm:tg}]
The proof is inspired by the proofs of \cite[Theorem~1]{McCormick2023}, \cite[Theorem~7.2]{McCormick2021}, and \cite[Theorem~4.5]{Tamstorf2021}.
Let $\vec{f}$, $\vec{v}^{[1]}$, $\vec{r}^{[1]}$, $\vec{r}^{[2]}_{\C}$, $\vec{v}^{[2]}_{\C}$, $\vec{v}^{[2]}$, $\vec{v}^{[3]}$, $\vec{v}^{[4]}$ denote the terms in \Cref{alg:two-grid} computed in exact arithmetic and let $\hat{\vec{f}}$, $\hat{\vec{v}}^{[1]}$, $\hat{\vec{r}}^{[1]}$, $\hat{\vec{r}}^{[2]}_{\C}$, $\hat{\vec{v}}^{[2]}_{\C}$, $\hat{\vec{v}}^{[2]}$, $\hat{\vec{v}}^{[3]}$, $\hat{\vec{v}}^{[3]}$ denote the corresponding terms computed in finite precision.
We first present a series of inequalities and bounds which are used below. They hold due to the assumptions imposed on the TG cycle:
\begin{align}
\label{eq:special_projection_A_Cnorm_bound}
\| \A^{-1}_\C \matrx{P}^{\top} \vec{v} \|_{\A_\C} & \leq  \|  \A^{-1}    \vec{v}\|_\A  \quad \forall \vec{v}, \\
\label{eq:B_C_Anorm_bound}
\| \matrx{M}_\C \A_{\C}  \|_{\A_\C}  & < 2, \\
\label{eq:f_bound}
 \| \vec{f}  \| & \leq \| \A \|^{\frac{1}{2}} \| \vec{y} \|_{\A}, \\
\label{eq:v[4]_Anorm_bound}
 \| \vec{v}^{[4]} \|_{\A}  & \leq 2 \| \vec{y}\|_\A ,\\
\label{eq:A-1_r^[1]_bound}
\| \A^{-1}  \vec{r}^{[1]}\|_{\A} & \leq \| \vec{y} \|_\A,
\\
\label{eq:r^[1]_bound}
\|  \vec{r}^{[1]}\| &\leq \| \A \|^{\frac{1}{2}} \| \vec{y}\|_\A,\\
\label{eq:A-1_Cr[2]_C_bound}
\| \A^{-1}_\C \vec{r}^{[2]}_\C \|_{\A_\C}  &\leq  \| \vec{y}\|_\A, \\
\label{eq:v[2]_C_A_C_norm_bound}
\| \vec{v}^{[2]}_{\C} \|_{\A_\C} &\leq 2 \| \vec{y} \|_\A.
\end{align}
Variants of the bounds \eqref{eq:special_projection_A_Cnorm_bound}-\eqref{eq:v[2]_C_A_C_norm_bound} can be found, e.g., in \cite{McCormick2023}. We include their derivation in \Cref{ape:mg} for self consistency of the text.

We focus on deriving the bound on the $\A$-norm of the error caused by computation in finite precision arithmetic in the TG scheme \eqref{eq:thm:tg_difference}. Analogous to the proof of \cite[Theorem~1]{McCormick2023}, we go line by line in \Cref{alg:two-grid} and bound the finite precision errors. Since some of the assumptions or bounds we use contain the Euclidean norm and some the $\A$-norm, we switch between these norms frequently in the derivation.

Line 1: Rounding $\vec{f}$ to $\dot{\varepsilon}$-precision arithmetic results in 
$\hat{\vec{f}} = \vec{f} + \delta_f$, where, using \eqref{eq:fp_round} and \eqref{eq:f_bound},
\begin{equation}\label{eq:delta_f_bound}
    \| \delta_f \| 
    \leq \dot{\varepsilon} \| \vec{f} \| 
    \leq \underbrace{\dot{\varepsilon} \| \A \|^{\frac{1}{2}}}_{K_0} \| \vec{y}  \|_{\A}.
\end{equation}

Line 2: Applying the smoothing to $\hat{\vec{f}} = \vec{f} + \delta_f$ in finite precision results in $\hat{\vec{v}}^{[1]} = \matrx{M}(\vec{f} + \delta_f) + \delta_{v^{[1]}}$, where, using assumption \eqref{eq:smoother_assumption_fp}, \eqref{eq:f_bound}, and \eqref{eq:delta_f_bound},
\begin{equation}\label{eq:delta_v[1]_bound}
    \| \delta_{v^{[1]}} \| \leq \Lambda_M (\| \vec{f}\| + \| \delta_f \|) \leq \underbrace{(\Lambda_M \| \A \|^\frac{1}{2}   + \Lambda_M K_0)}_{K_1}  \| \vec{y}  \|_{\A}.
\end{equation}
The term  $\hat{\vec{v}}^{[1]}$ can be written as $ \hat{\vec{v}}^{[1]}=\vec{v}^{[1]} + \Delta_{v^{[1]}}$, where $\Delta_{v^{[1]}}= \matrx{M} \delta_f +  \delta_{v^{[1]}}$ is the accumulated error and, using \eqref{eq:delta_f_bound} and \eqref{eq:delta_v[1]_bound},
\begin{equation}\label{eq:Delta_v[1]_bound}
    \| \Delta_{v^{[1]}} \| \leq \| \matrx{M} \delta_f \| + \|  \delta_{v^{[1]}} \| \leq \underbrace{(\| \matrx{M} \| K_0 + K_1)}_{K_2}  \| \vec{y}  \|_{\A}.
\end{equation}

Line 3: Computing 
$(\vec{f} + \delta_f) - \A (\vec{v}^{[1]}+\Delta_{v^{[1]}})$
in $\dot{\varepsilon}$-precision results in 
$\hat{\vec{r}}^{[1]}=\vec{f} + \delta_f - \A(\vec{v}^{[1]}+\Delta_{v^{[1]}}) + \delta_{r^{[1]}}$,    
where, using \eqref{eq:fp_mat_vec_add_quant},  $\vec{v}^{[1]}=\matrx{M}\vec{f}$,  \eqref{eq:f_bound}, \eqref{eq:delta_f_bound}, and \eqref{eq:Delta_v[1]_bound},
\begin{align}
\nonumber
\| \delta_{r^{[1]}} \|
&\leq  \dot{\varepsilon} m_{A,\dot{\varepsilon}}
( \| \vec{f} + \delta_f \| + \| | \A | \|  \cdot \|  \vec{v}^{[1]}+\Delta_{v^{[1]}}\| ) 
 \\ \nonumber & 
\leq  \dot{\varepsilon} m_{A,\dot{\varepsilon}}
( \| \vec{f}  \| + \| \delta_f \| + \| | \A | \|( \|  \vec{v}^{[1]} \| + \| \Delta_{v^{[1]}}\| ) ) \\
\label{eq:delta_r[1]_bound}
& \leq \underbrace{ \dot{\varepsilon} m_{A,\dot{\varepsilon}}
( \| \A \|^{\frac{1}{2}} + K_0 + \| | \A | \|( \| \matrx{M} \|\| \A\|^{\frac{1}{2}} + K_2  ))}_{K_3} \| \vec{y}  \|_{\A}. 
\end{align}
The term  $\hat{\vec{r}}^{[1]}$ can be written as $ \hat{\vec{r}}^{[1]}=\vec{r}^{[1]}+ \Delta_{r^{[1]}}$, where 
$\Delta_{r^{[1]}} = \delta_f - \A \Delta_{v^{[1]}} + \delta_{r^{[1]}}$
is the accumulated error, and using \eqref{eq:delta_f_bound}, \eqref{eq:Delta_v[1]_bound}, and \eqref{eq:delta_r[1]_bound},
\begin{align}
\nonumber
\| \Delta_{r^{[1]}} \| &=  \| \delta_f - \A \Delta_{v^{[1]}} + \delta_{r^{[1]}} \| 
\leq \| \delta_f \| +  \| \A \|   \| \Delta_{v^{[1]}}\|  + \| \delta_{r^{[1]}} \| \\
\label{eq:Delta_r[1]_bound}
&\leq \underbrace{(K_0 + \| \A \|   K_2 + K_3 )}_{K_4} \| \vec{y} \|_\A.
\end{align}

Line 4: Computing $ \matrx{P}^{\top}(\vec{r}^{[1]} + \Delta_{r^{[1]}})$ in $\dot{\varepsilon}$-precision results in 
$\hat{\vec{r}}^{[2]}_\C=\matrx{P}^{\top}(\vec{r}^{[1]} + \Delta_{r^{[1]}}) + \delta_{r^{[2]}_{\C}}$, 
where, using \eqref{eq:fp_mat_vec_quant}, \eqref{eq:r^[1]_bound}, and \eqref{eq:Delta_r[1]_bound},
\begin{equation}
\| \delta_{r^{[2]}_{\C}} \| 
 \leq 
\dot{\varepsilon} \underline{m}_{P,\dot{\varepsilon}}\| | \matrx{P} | \| ( \| \vec{r}^{[1]} \|  +  \| \Delta_{r^{[1]}} \| )
\label{eq:delta_r[2]_C_bound}
 \leq 
\underbrace{\dot{\varepsilon} \underline{m}_{P,\dot{\varepsilon}}\| | \matrx{P} | \| ( \| \A \|^{\frac{1}{2}}   
+  K_4 ) }_{K_5}\| \vec{y}  \|_{\A}.
\end{equation}
The term $\hat{\vec{r}}^{[2]}_\C$ can be written as $\hat{\vec{r}}^{[2]}_\C=\vec{r}^{[2]}_\C + \Delta_{r^{[2]}_\C}$, where $\Delta_{r^{[2]}_\C} = \matrx{P}^{\top} \Delta_{r^{[1]}} + \delta_{r^{[2]}_{\C}}$ is the accumulated error and, using \eqref{eq:special_projection_A_Cnorm_bound},  \eqref{eq:A^-1_v_by_Euclid}, $\| \matrx{P}^{\top} \|=\| \matrx{P}\| $, 
\eqref{eq:Delta_v[1]_bound}, \eqref{eq:delta_f_bound},  \eqref{eq:delta_r[1]_bound}, and \eqref{eq:delta_r[2]_C_bound},
\begin{align}
\nonumber
\|\A^{-1}_{\C} \Delta_{r^{[2]}_\C}\|_{\A_\C} & =
\| \A^{-1}_{\C}(\matrx{P}^{\top} \Delta_{r^{[1]}}  +  \delta_{r^{[2]}_{\C}}) \|_{\A_{\C}} \\
\nonumber
&= \| \A^{-1}_{\C}(\matrx{P}^{\top}( \delta_f - \A \Delta_{v^{[1]}} + \delta_{r^{[1]}})  +  \delta_{r^{[2]}_{\C}}) \|_{\A_{\C}} \\
\nonumber
& \leq  \| \A^{-1}_{\C}\matrx{P}^{\top} \A \Delta_{v^{[1]}} \|_{\A_{\C}} 
+ \| \A^{-1}_{\C}(\matrx{P}^{\top}( \delta_f + \delta_{r^{[1]}})  +  \delta_{r^{[2]}_{\C}}) \|_{\A_{\C}} \\
\nonumber
&\leq \| \Delta_{v^{[1]}} \|_{\A} + 
\| \A^{-1}_{\C} \|^{\frac{1}{2}} \| \matrx{P}^{\top} (  \delta_f  +  \delta_{r^{[1]}} ) + \delta_{r^{[2]}_{\C}} \| \\
\nonumber
&\leq \|  \Delta_{v^{[1]}} \|_{\A} + 
\| \A^{-1}_{\C} \|^{\frac{1}{2}} (\| \matrx{P}\| ( \| \delta_f \| + \| \delta_{r^{[1]}}\| )  + \| \delta_{r^{[2]}_{\C}} \|) \\
\label{eq:Delta_r^[2]_C_bound}
&\leq  \underbrace{ ( \| \A \|^{\frac{1}{2}} K_2 + \| \A^{-1}_{\C} \|^{\frac{1}{2}} (\| \matrx{P}\| ( K_0 + K_3 )  + K_5 )) }_{K_6} \| \vec{y}  \|_{\A}.
\end{align}

Line 5: Applying the coarse-grid solver to $\vec{r}^{[2]}_\C + \Delta_{r^{[2]}_{\C}}$ in finite precision results in
$\hat{\vec{v}}^{[2]}_{\C} = \matrx{M}_{\C} (\vec{r}^{[2]}_\C + \Delta_{r^{[2]}_\C}) + \delta_{v^{[2]}_\C}$,
where, using the assumption \eqref{eq:solver_assumption_fp} and the estimates \eqref{eq:A-1_Cr[2]_C_bound} and \eqref{eq:Delta_r^[2]_C_bound},
\begin{align}
\nonumber
\| \delta_{v^{[2]}_\C} \|_{\A_\C} & 
\leq \Lambda_{\C} \| \A^{-1}_{\C} ( \vec{r}^{[2]}_\C + \Delta_{r^{[2]}_\C} ) \|_{\A_\C} 
\leq \Lambda_\C ( \| \A^{-1}_{\C} \vec{r}^{[2]}_\C \|_{\A_\C} +  \|\A^{-1}_{\C} \Delta_{r^{[2]}_\C}\|_{\A_\C}) \\
\label{eq:delta_v[2]_C_bound}
&\leq \underbrace{\Lambda_\C ( 1 +  K_6 )}_{K_{7}}\| \vec{y}\|_{\A}.
\end{align}
The term $\hat{\vec{v}}^{[2]}_{\C}$ can be written as
$\hat{\vec{v}}^{[2]}_{\C}=\vec{v}^{[2]}_{\C} + \Delta_{v^{[2]}_{\C}}$, 
where 
$\Delta_{v^{[2]}_{\C}} = \matrx{M}_{\C} \Delta_{r^{[2]}_{\C}} + \delta_{v^{[2]}_\C}$ is the  accumulated error and, using \eqref{eq:B_C_Anorm_bound}, \eqref{eq:Delta_r^[2]_C_bound}, and \eqref{eq:delta_v[2]_C_bound},
\begin{align}
\nonumber
\| \Delta_{v^{[2]}_{\C}} \|_{\A_{\C}}
& \leq  \| \matrx{M}_{\C} \matrx{A}_{\C}\matrx{A}^{-1}_{\C}  \Delta_{r^{[2]}_{\C}} \|_{\A_{\C}} +   \| \delta_{v^{[2]}_\C} \|_{\A_{\C}} 
 \\ \nonumber & 
\leq \| \matrx{M}_{\C} \A_\C  \|_{\A_{\C}} 
\| \A^{-1}_{\C} \Delta_{r^{[2]}_{\C}}  \|_{\A_{\C}} + \| \delta_{v^{[2]}_\C} \|_{\A_{\C}} \\
\label{eq:Delta_v[2]_C}
&\leq  \underbrace{ (2 K_6 + K_7)}_{K_8} \| \vec{y}  \|_{\A}.
\end{align}

Line 6: 
Computing $ \matrx{P} (\vec{v}^{[2]}_{\C} + \Delta_{v^{[2]}_{\C}})$ in $\dot{\varepsilon}$-precision results in 
$\hat{\vec{v}}^{[3]}=\matrx{P}(\vec{v}^{[2]}_{\C} + \Delta_{v^{[2]}_{\C}}) +  \delta_{v^{[3]}}$, 
where, using \eqref{eq:A_norm_by_Euclid}, \eqref{eq:fp_mat_vec_quant}, \eqref{eq:Euclid_by_A_norm}, $\kappa^{\frac{1}{2}}_A=\| \A \|^{\frac{1}{2}} \| \A^{-1} \|^{\frac{1}{2}}$, $\xi = \| \A \|^{\frac{1}{2}}/ \| \A_\C \|^{\frac{1}{2}}$,  \eqref{eq:v[2]_C_A_C_norm_bound}, and \eqref{eq:Delta_v[2]_C},
\begin{align}
\nonumber
\| \delta_{v^{[3]}} \|_{\A}  &\leq \| \A \|^{\frac{1}{2}}  \| \delta_{v^{[3]}} \| 
\leq \| \A \|^{\frac{1}{2}}\dot{\varepsilon} \underline{m}_{P,\dot{\varepsilon}}
\|| \matrx{P}| \|  \| \vec{v}^{[2]}_{\C} + \Delta_{v^{[2]}_{\C}} \|  \\
\nonumber
& \leq \| \A \|^{\frac{1}{2}} 
\dot{\varepsilon} \underline{m}_{P,\dot{\varepsilon}}
\|| \matrx{P}| \| 
\|  \A^{-1}_{\C} \|^{\frac{1}{2}}
\| \vec{v}^{[2]}_{\C} + \Delta_{v^{[2]}_{\C}} \|_{\A_{\C}}\\
\nonumber
&\leq \dot{\varepsilon}\kappa^{\frac{1}{2}}_{A} \xi  \underline{m}_{P,\dot{\varepsilon}}  
\|| \matrx{P}| \| 
(\| \vec{v}^{[2]}_{\C} \|_{\A_{\C}} + \| \Delta_{v^{[2]}_{\C}} \|_{\A_{\C}})\\
\nonumber
& \leq \dot{\varepsilon}\kappa^{\frac{1}{2}}_{A} \xi  \underline{m}_{P,\dot{\varepsilon}}  
\|| \matrx{P}| \| 
(2 \| \vec{y} \|_{\A} + \| \Delta_{v^{[2]}_{\C}} \|_{\A_{\C}})\\
\label{eq:delta_v[3]_bound}
& \leq \underbrace{\dot{\varepsilon}\kappa^{\frac{1}{2}}_{A} \xi  \underline{m}_{P,\dot{\varepsilon}}  
\|| \matrx{P}| \| 
(2 + K_8)}_{K_9}\| \vec{y} \|_{\A}. 
\end{align}
The term $\hat{\vec{v}}^{[3]}$ can be written as $\hat{\vec{v}}^{[3]}=\vec{v}^{[3]} + \Delta_{v^{[3]}}$, where 
$\Delta_{v^{[3]}}=\matrx{P} \Delta_{v^{[2]}_{\C}} +\delta_{v^{[3]}} $
is the accumulated error, and using $\A_{\C} = \matrx{P}^{\top}\A \matrx{P}$, \eqref{eq:Delta_v[2]_C}, and \eqref{eq:delta_v[3]_bound},
\begin{align}
\nonumber
\|  \Delta_{v^{[3]}} \|_\A 
&= \| \matrx{P} \Delta_{v^{[2]}_{\C}} +\delta_{v^{[3]}} \|_\A 
\leq \| \matrx{P}\Delta_{v^{[2]}_{\C}} \|_{\A} + \| \delta_{v^{[3]}} \|_\A 
= \| \Delta_{v^{[2]}_{\C}} \|_{\A_{\C}} + \| \delta_{v^{[3]}} \|_\A \\
\label{eq:Delta_v[3]_bound}
&\leq \underbrace{(K_8 + K_9)}_{K_{10}} \| \vec{y} \|_{\A}. 
\end{align}

Line 7: Computing $\vec{v}^{[1]} +\Delta_{v^{[1]}} + \vec{v}^{[3]}  + \Delta_{v^{[3]}}$ in $\dot{\varepsilon}$-precision results in 
$ \hat{\vec{v}}^{[4]} = \vec{v}^{[1]} +\Delta_{v^{[1]}} + \vec{v}^{[3]}  + \Delta_{v^{[3]}} + \delta_{v^{[4]}},$
where, using  \eqref{eq:A_norm_by_Euclid}, \eqref{eq:fp_add}, $  \vec{v}^{[4]}=  \vec{v}^{[1]} + \vec{v}^{[3]}$, 
\eqref{eq:Euclid_by_A_norm}, 
\eqref{eq:v[4]_Anorm_bound}, \eqref{eq:Delta_v[1]_bound}, and \eqref{eq:Delta_v[3]_bound},
\begin{align}
\nonumber
\| \delta_{v^{[4]}} \|_\A & \leq \| \A \|^{\frac{1}{2}} \| \delta_{v^{[4]}} \|  \\
\nonumber
& \leq  \| \A \|^{\frac{1}{2}} \dot{\varepsilon} \| \vec{v}^{[1]} +\Delta_{v^{[1]}} + \vec{v}^{[3]}  + \Delta_{v^{[3]}} \| \\
\nonumber
&\leq  \| \A \|^{\frac{1}{2}} \dot{\varepsilon} (\| \vec{v}^{[4]} \| + \| \Delta_{v^{[1]}} \|   +\| \Delta_{v^{[3]}} \|) \\
\nonumber
& \leq  \| \A \|^{\frac{1}{2}} \dot{\varepsilon} ( \| \A^{-1} \|^{\frac{1}{2}} \| \vec{v}^{[4]} \|_{\A} + \| \Delta_{v^{[1]}} \|  + \| \A^{-1} \|^{\frac{1}{2}} \| \Delta_{v^{[3]}} \|_{\A}) \\
\label{eq:delta_v[4]_bound}
& \leq \underbrace{\dot{\varepsilon} ( 2 \kappa^{\frac{1}{2}}_A  +  \| \A \|^{\frac{1}{2}} K_2  + \kappa^{\frac{1}{2}}_A K_{10})}_{K_{11}} \| \vec{y} \|_{\A}. 
\end{align}
Finally the computed approximation $\hat{\vec{v}}^{[4]}$ can be written as $\hat{\vec{v}}^{[4]}=\vec{v}^{[4]} + \Delta_{v^{[4]}}$, where
$\Delta_{v^{[4]}} = \Delta_{v^{[1]}} + \Delta_{v^{[3]}} + \delta_{v^{[4]}}$ 
is the accumulated error and using \eqref{eq:Delta_v[1]_bound}, \eqref{eq:Delta_v[3]_bound}, \eqref{eq:A_norm_by_Euclid}, and \eqref{eq:delta_v[4]_bound},
\begin{align*}
     \| \Delta_{v^{[4]}}\|_{\A} & \leq  \| \Delta_{v^{[1]}} \|_{\A} + \| \Delta_{v^{[3]}} \|_{\A} + \| \delta_{v^{[4]}} \|_{\A}
     \leq \underbrace{(\| \A \|^{\frac{1}{2}} K_2 + K_{10} + K_{11} )}_{\lambda_{\TG}} \| \vec{y} \|_{\A}.
\end{align*}
Since $\hat{\vec{v}}^{[4]}=\hat{\vec{y}}_{\TG}$ and $\vec{v}^{[4]}=\vec{y}_{\TG}$, we have $\Delta_{v^{[4]}} = \hat{\vec{y}}_{\TG} - \vec{y}_{\TG} = \delta_{\TG}$. Then
\begin{equation*}
\| \vec{y}_{\TG} - \hat{\vec{y}}_{\TG}  \|_\A  \leq \lambda_{\TG} \| \vec{y} \|_\A. 
\end{equation*}
We rewrite the expression for $\lambda_{\TG}$ by grouping higher order terms in a remainder~$R$. The higher order terms are negligible in a local sense and include, e.g., second powers of $\dot{\varepsilon}$, $\Lambda_\C$, $\Lambda_M$, or a product of at least two of them. 
All remainders $R_k$, $k=1,\ldots,6$, defined below contain only high order terms.
Listing and rewriting the constants $K_j$, $j=0,\ldots,11$, leads to
\begin{flalign*}
K_0& = \dot{\varepsilon} \| \A \|^{\frac{1}{2}}, &&& \\
K_1& = (\Lambda_M \| \A \|^\frac{1}{2}   + \Lambda_M K_0) = (\Lambda_M  + \Lambda_M\dot{\varepsilon}) \| \A \|^{\frac{1}{2}}, &&& \\
K_2& = \| \matrx{M} \| K_0 + K_1 = (\| \matrx{M} \| \dot{\varepsilon} + \Lambda_M +\Lambda_M\dot{\varepsilon}  )\| \A \|^{\frac{1}{2}},  &&& \\ 
K_3& =  \dot{\varepsilon} m_{A,\dot{\varepsilon}} ( \| \A \|^{\frac{1}{2}} + K_0 + \| | \A | \|( \| \matrx{M} \|\| \A\|^{\frac{1}{2}} + K_2  )) &&& \\
& = \dot{\varepsilon} m_{A,\dot{\varepsilon}} \| \A \|^{\frac{1}{2}}( 1 + \| | \A | \| \| \matrx{M} \|) + \underbrace{\dot{\varepsilon} m_{A,\dot{\varepsilon}} (K_0 + \| | \A | \|K_2)}_{R_1},
\end{flalign*}
\begin{flalign*}
K_4 & = K_0 + \| \A \|   K_2 + K_3  &&& \\
& = (\dot{\varepsilon} + \| \A \| ( \| \matrx{M} \| \dot{\varepsilon} + \Lambda_M + \Lambda_M \dot{\varepsilon} ))\| \A \|^{\frac{1}{2}} + \dot{\varepsilon}m_{A,\dot{\varepsilon}} \| \A\|^{\frac{1}{2}}( 1 + \| | \A | \| \| \matrx{M} \|   )  + R_1, &&& \\
K_5 & = \dot{\varepsilon} \underline{m}_{P,\dot{\varepsilon}} \| | \matrx{P} | \| ( \| \A \|^{\frac{1}{2}}   
+  K_4 )
 = \dot{\varepsilon} \underline{m}_{P,\dot{\varepsilon}}  \| | \matrx{P} | \| \| \A \|^{\frac{1}{2}} + 
\underbrace{ \dot{\varepsilon} \underline{m}_{P,\dot{\varepsilon}}  \| | \matrx{P} | \| K_4 }_{R_{2}}, &&&  \\
K_6 & =   \| \A \|^{\frac{1}{2}} K_2 + \| \A^{-1}_{\C} \|^{\frac{1}{2}} (\| \matrx{P}\| ( K_0 + K_3 )  + K_5 ) &&& \\
& = \| \A \|  (\| \matrx{M} \| \dot{\varepsilon} + \Lambda_M ) + \dot{\varepsilon} \kappa^{\frac{1}{2}}_{A} \xi (\| \matrx{P}\| (1+m_{A,\dot{\varepsilon}}(1+  \| | \A | \| \| \matrx{M} \|)) + \underline{m}_{P,\dot{\varepsilon}} \| | \matrx{P} | \| ) &&& \\
& \qquad 
+ \underbrace{\| \A \|  \Lambda_M \dot{\varepsilon}  + \| \A^{-1}_{\C} \|^{\frac{1}{2}}  \| \matrx{P} \|  R_1 + \| \A^{-1}_{\C} \|^{\frac{1}{2}}  R_2}_{R_3}, &&& \\
K_7 & = \Lambda_\C ( 1 +  K_6 ) = \Lambda_\C + \underbrace{\Lambda_\C K_6}_{R_4}, &&& 
\end{flalign*}
\begin{flalign*}
K_8& = 2K_6 + K_7  = 2K_6 + \Lambda_\C + R_4, &&& \\
K_9 & = \dot{\varepsilon}\kappa^{\frac{1}{2}}_{A} \xi  \underline{m}_{P,\dot{\varepsilon}} 
\|| \matrx{P}| \|
(2 + K_8)
= 2\dot{\varepsilon}\kappa^{\frac{1}{2}}_{A} \xi  \underline{m}_{P,\dot{\varepsilon}} 
\|| \matrx{P}| \| +
\underbrace{ \dot{\varepsilon} \kappa^{\frac{1}{2}}_{A} \xi \underline{m}_{P,\dot{\varepsilon}} \|| \matrx{P}| \|K_8  }_{R_5}, &&& \\
K_{10} & = K_8 +K_9 
&&& \\
K_{11}&= \dot{\varepsilon} (2\kappa^{\frac{1}{2}}_A  +  \| \A \|^{\frac{1}{2}} K_2  + \kappa^{\frac{1}{2}}_A K_{10})
= 2 \dot{\varepsilon} \kappa^{\frac{1}{2}}_A  +  \underbrace{\dot{\varepsilon}\| \A \|^{\frac{1}{2}} K_2  + \dot{\varepsilon}\kappa^{\frac{1}{2}}_A K_{10}}_{R_6}.    &&&
\end{flalign*}
Finally,  $\lambda_{\TG}$ can be rewritten as 
\begin{align*}
\lambda_{\TG}& = \| \A \|^{\frac{1}{2}} K_2 + K_{10} +  K_{11} \\
& = \| \A \| (\| \matrx{M} \| \dot{\varepsilon} + \Lambda_M ) + 
2K_6 + \Lambda_\C + R_4 + 
2\dot{\varepsilon} \kappa^{\frac{1}{2}}_{A} \xi\underline{m}_{P,\dot{\varepsilon}}  
 \|| \matrx{P}| \| 
 + R_5
+ 2\dot{\varepsilon}\kappa^{\frac{1}{2}}_{A} +  R_6\\
& = 3\| \A \| (\| \matrx{M} \| \dot{\varepsilon} + \Lambda_M  ) + 
\dot{\varepsilon} \kappa^{\frac{1}{2}}_{A} \xi (2\| \matrx{P}\| (1+m_{A,\dot{\varepsilon}})(1+  \| | \A | \| \| \matrx{M} \|)) 
+ 4\underline{m}_{P,\dot{\varepsilon}} \| | \matrx{P} | \| ) \\
& \qquad \qquad  + \Lambda_\C + 2\dot{\varepsilon}\kappa^{\frac{1}{2}}_{A} + \underbrace{2R_3 + R_4  + R_5 + R_6}_{R}.
\end{align*}
Since $ \kappa^{\frac{1}{2}}_{A} \geq 1 $ the term $3\| \A \| \| \matrx{M} \| \dot{\varepsilon}$ can be bounded by $3\| \A \| \| \matrx{M} \| \dot{\varepsilon} \kappa^{\frac{1}{2}}_{A}$.
Defining $C_1$ and $C_2$ as
\begin{equation*}
C_1  = 2\xi \| \matrx{P} \| (1 + m_{A,\dot{\varepsilon}}) \| | \matrx{A} | \|   + 3 \|  \matrx{A}  \|, \quad
C_2  = \xi ( 2\| \matrx{P} \| (1 + m_{A,\dot{\varepsilon}}) +  4\underline{m}_{P,\dot{\varepsilon}} \| | \matrx{P} | \|) + 2,
\end{equation*}
gives \eqref{eq:thm:tg_difference}.
Using the assumption \eqref{eq:TG_exact_rho} and \eqref{eq:thm:tg_difference} we have
\begin{equation*}      
 \| \vec{y} - \hat{\vec{y}}_{\TG} \|_\A \leq  \| \vec{y} - \vec{y}_{\TG} \|_\A + \| \vec{y}_{\TG} - \hat{\vec{y}}_{\TG}  \|_\A \leq (\rho_\TG + \Lambda_\TG) \| \vec{y} \|_\A,
\end{equation*}
which finishes the proof.
\end{proof}

\section{V-cycle method}
\label{sec:V-cycle}
In this section, we present the finite precision analysis of a V-cycle method which can be seen as a generalization of the TG method.

We consider using the V-cycle method for solving  $\A\vec{y}=\vec{f}$ and present a bound on the finite precision error of one V-cycle starting with zero initial approximation. As before, the analysis of multiple cycles can be obtain by viewing the method as a iterative refinement method.

In the V-cycle method, the approximate solution is computed using a hierarchy of $J+1$ levels numbered from $0$ to $J$ from the coarsest to the finest level. Each level contains a system matrix $\A_j \in \R^{n_j \times n_j}$, $j=0,\ldots,J$, with $\A_J=\A$.
The information is transferred between the $(j-1)$th and $j$th levels using a full-rank prolongation matrix $\matrx{P}_j\in \R^{n_j \times n_{j-1}}$, $j=1,\ldots,J$. The restriction matrices are transposes of the prolongation matrices. We assume that the system matrices satisfy the Galerkin condition,  i.e., $\A_{j-1} = \matrx{P}^{\top}_j \A_j \matrx{P}_j$, $j=1,\ldots,J$.

The computation consists of smoothing on fine levels and  solving a system of linear equations on the coarsest-level.
Analogously, as in the TG cycle, we consider, for simplicity of the analysis, that the smoothing is applied before the computation on the coarser levels. We assume that the application of the coarsest level solver or the smoother on level $j$, $j=0,1,\ldots,J$, in exact arithmetic can, for any vector $\vec{f}_j \in \R^{n_j}$, be expressed as $\matrx{M}_j \vec{f}_j$, where $\matrx{M}_j \in \R^{n_j \times n_j}$ is a non-singular matrix and
\begin{equation}\label{eq:V_cycle_smoothing_A_norm_assumption}
\| \matrx{I}_j - \matrx{M}_j \A_j \|_{\A_j}<1.
\end{equation}

We assume that all operations on a fine level $j$, $j=1,\dots,J$, besides the smoothing, are done in finite precision arithmetic with unit a roundoff $\dot{\varepsilon}_j$. 
We assume that the precision used on level $j$, $j=2,\ldots,J$, is higher or equal to the precision used on the coarser level $j-1$, i.e., $\dot{\varepsilon}_j \leq \dot{\varepsilon}_{j-1}$, and that $\dot{\varepsilon}_j$ is small enough such that $(m_{A_j}+2) \dot{\varepsilon}_j < 1$ and $(\underline{m}_{P_j}+1) \dot{\varepsilon}_j < 1$.
We assume that the matrices $\A_j$ and $\matrx{P}_j$ on level $j$ are rounded to  $\dot{\varepsilon}_j$-precision for computing the residual, restriction, and prolongation.

We assume that there exists a positive constant $\Lambda_{M_j}$, $j=1\ldots,J$, and $\Lambda_0$ such that for any vector $\vec{f}_j$, the application of smoothing in finite precision on the $j$th level and the application of the coarsest level solver results in, respectively,
\begin{align}
\label{eq:smoother_assumption_fp_V-cycle}
& \matrx{M}_j \vec{f}_j + \delta_{M_j}, \quad \| \delta_{M_j} \| \leq \Lambda_{M_j}  \| \vec{f}_j\|, \\
\label{eq:solver_assumption_fp_v_cycle}
& \matrx{M}_0 \vec{f}_0 + \delta_0, \quad \| \delta_{0} \|_{\A_0} \leq \Lambda_{0} \| \A^{-1}_{0} \vec{f}_0\|_{\A_{0}},
\end{align}
and the results belong to the $\dot{\varepsilon}_j$-precision arithmetic and $\dot{\varepsilon}_1$-precision arithmetic for the coarsest-level solver.

A recursive formulation of the V-cycle starting with zero initial approximation is described in \Cref{alg:V-cycle}. 
We assume that there exists a factor $\rho_\V<1$, such that for any $j=1,\ldots,J$, the V-cycle in exact arithmetic with $j+1$ levels $0,\ldots,j$ reduces the $\A_j$-norm of the error at least by the factor $\rho_\V$. In more detail, let $\vec{f}_j$, $j=1,\ldots,J$, be a right-hand side vector on the $j$th level, let $\vec{y}_j$, $j=1,\ldots,J$, be the solution of $\A_j\vec{y}_j = \vec{f}_j$, and let $\vec{y}_{\V,j}$ be the approximation computed using the V-cycle  in exact arithmetic with $j+1$ levels, $0,\ldots,j$. Then 
\begin{equation}\label{eq:uniform_convergence_V_cycle}
    \| \vec{y}_j - \vec{y}_{\V,j}  \|_{\A_j} \leq \rho_{\V} \| \vec{y}_j \|_{\A_j}.
\end{equation}

\begin{algorithm}[H]
\caption{V-cycle with zero initial approximation, $\mathbf{V}( \vec{f}_{j},j$).}
\label{alg:V-cycle}
\begin{algorithmic}[1]
\IF{$ j \neq 0 $}
\STATE{$\vec{f}_j \leftarrow\mathrm{Round}(\vec{f}_j,\vp\text{-precision})$ }
\STATE{$\vec{v}^{[1]}_j \leftarrow \matrx{M}_j \vec{f}_j$ \COMMENT{Apply smoothing.}}
\STATE{$\vec{r}^{[1]}_j \leftarrow \vec{f}_j -  \A_j \vec{v}^{[1]}_j$ \COMMENT{Compute residual in $\vp$-precision.}}
\STATE{$\vec{r}^{[2]}_{j-1} \leftarrow \matrx{P}_j^{\top} \vec{r}^{[1]}_j$ \COMMENT{Restrict the residual to level $j-1$ in $\vp$-precision.}}
\STATE{$\vec{v}^{[2]}_{j-1} \leftarrow \matrx{V}( \vec{r}^{[2]}_{j-1},
j-1)$ \COMMENT{Recursive call.}}
\STATE{$\vec{v}^{[3]}_j \leftarrow \matrx{P}_j  \vec{v}^{[2]}_{j-1}$ \COMMENT{Prolongate the correction to level $j$ in $\vp$-precision.}}
\STATE{$\vec{v}^{[4]}_j \leftarrow \vec{v}^{[1]}_j + \vec{v}^{[3]}_j$ \COMMENT{Correct the previous approximation in $\vp$-precision.}}
\RETURN{$\vec{v}^{[4]}_j$}
\ELSE
\RETURN{$\matrx{M}_0 \vec{f}_0$ \COMMENT{Coarsest-level solve.
}}
\ENDIF
\end{algorithmic}
\end{algorithm}

We present the following result on the effects of finite precision errors in the
V-cycle. Its proof, based on consecutive usage of \Cref{thm:tg}, is presented below.
\begin{theorem}\label{thm:V-cycle}
Let $\vec{y}_{\V}$ and  $\hat{\vec{y}}_{\V}$ be the approximate solution of $\A\vec{y}=\vec{f}$ computed using one V-cycle (\Cref{alg:V-cycle}) with $J+1$ levels, applied exactly and in finite precision, respectively.
The $\A$-norm of the finite precision error $\delta_{\V}=\vec{y}_{\V} - \hat{\vec{y}}_{\V}$ can be bounded as  $\|  \delta_{\V} \|_{\A} \leq \Lambda_{\V}  \| \vec{y} \|_{\A}$ with 
\begin{equation*}
\Lambda_{\V} = \Lambda_0 + \sum^{J} _{j=1} \Lambda_{V,j},  \quad \Lambda_{V,j} =    3\|  \A_j  \| \Lambda_{M_j} 
+  \dot{\varepsilon}_j \kappa^{\frac{1}{2}}_{A_j} (C_{1,j}   \| \matrx{M}_j \|
+ C_{2,j} )
+ R_j,
\end{equation*}
with $C_{1,j}$ and $C_{2,j}$ being positive constants depending on  $\| \A_j \|$, $\| | \A_j | \|$, $\| \matrx{P}_j\|$, $ \| | \matrx{P}_j | \|$, 
\begin{equation*}
m_{A_j, \varepsilon_j} = \frac{(m_{A_j}+2)}{1-(m_{A_j}+2)\dot{\varepsilon}_j}, \quad 
    \underline{m}_{P_j, \varepsilon_j} = \frac{(\underline{m}_{P_j}+1)}{1-(\underline{m}_{P_j}+1)\dot{\varepsilon}_j}, \quad \text{and} \quad 
  \frac{\| \A^{-1}_{j-1}\|^{\frac{1}{2}}} { \| \A^{-1}_j \|^{\frac{1}{2}}}.
\end{equation*}
The remainders $R_j$ contain additional higher order terms. If the sum of the terms in $\Lambda_{\V}$ excluding the quantities $R_j$ is sufficiently smaller than one, the quantities $R_j$ are negligible. 
If~$\rho_{\V} + \Lambda_{\V} < 1$, the V-cycle applied in finite precision reduces the relative $\A$-norm of the error by at least the factor $ \rho_{\V} + \Lambda_{\V} $, i.e.,
$\| \vec{y} - \hat{\vec{y}}_{\V}  \|_{\A}  \leq (\rho_{\V} +  \Lambda_{\V}) \| \vec{y} \|_{\A}$.
In other words, the V-cycle applied in finite precision is a contraction with a contraction factor $ \rho_{\V} + \Lambda_{\V}$.
\end{theorem}
This theorem provides insight into how the finite precision errors coming from the coarsest-level solver, the smoothers, and from computing the residual, restriction, prolongation, and correction on  individual levels may affect the overall finite precision error. We see that the requirement on the $\dot{\varepsilon}_j$-precision as well as the finite precision error of the smoother may differ on each fine level based on the properties of the corresponding system, prolongation matrix, and the smoother.

\begin{proof}[Proof of \Cref{thm:V-cycle}]
We use induction on the number of levels.
The V-cycle with two levels can be seen as the TG cycle with $\matrx{M} = \matrx{M}_{1}$  and $\matrx{M}_\C = \matrx{M}_{0}$. Since the assumptions of \Cref{thm:tg} are satisfied, the statement holds for $j=1$.

Let $\matrx{V}_j$, $j=1,\ldots,J$, be the matrix corresponding to applying the V-cycle scheme in exact arithmetic with $j+1$ levels $0,\ldots,j$. Such a matrix exists; see, e.g., \cite[Theorem~2.4.1]{Trottenberg2001}. The assumption \eqref{eq:uniform_convergence_V_cycle} yields
\begin{equation}\label{eq:v-cycle_error_matrix_A_norm_bound}
\| \matrx{I}_j - \matrx{V}_j\A_j \|_{\A_j} = \max_{\vec{y}_j} \frac{ \| (\matrx{I}_j - \matrx{V}_j\A_j)\vec{y}_{j} \|_{\A_j}}{\| \vec{y}_j \|_{\A_j}} \leq \rho_V<1.
\end{equation}

We assume that the statement of the theorem holds for the V-cycle scheme with $j$ levels.
We can view the V-cycle scheme with $j+1$ levels as a two-grid correction scheme where the coarse-grid solver is the V-cycle scheme with $j$ levels, i.e., $\matrx{M}=\matrx{M}_j$ and $\matrx{M}_\C = \matrx{V}_{j-1}$. Since the smoothing routine on level $j$ and the coarse-grid solver satisfy the assumptions of \Cref{thm:tg}, in particular,
\begin{equation*}
\| \matrx{I}_\C-\matrx{M}_\C \A_{\C}\|_{\A_\C} = \| \matrx{I}_{j-1}-\matrx{V}_{j-1} \A_{j-1}\|_{\A_{j-1}} < 1, \qquad
    \Lambda_\C  = \Lambda_0 + \sum^{j-1}_{i=0}   \Lambda_{\V,j},
\end{equation*}
the result also holds for the V-cycle scheme with $j+1$ levels.
\end{proof}

\section{Smoothing based on incomplete Cholesky factorization}
\label{sec:IC_smoothing}

In this section, we present a mixed precision smoothing routine based on incomplete Cholesky (IC) factorization and a finite precision error analysis of its application. 
These results are combined with the results from previous section in \Cref{sec:IC-V-cycle}, where we summarize the results on the effects of finite precision errors in the V-cycle with IC smoothing.

We assume that an SPD matrix $\A$ is approximated by its IC factorization $\matrx{L}  \matrx{L}^{\top}$, with $\matrx{L}$ being a lower triangular matrix; see, e.g., \cite[Chapter~10]{Saad2003}, \cite[Chapter~10]{Scott2023}.
The application of the smoother involves solving triangular systems with the matrix $\matrx{L}$ and its transpose, which we assume is done using substitution; see, e.g., \cite[Chapter~8]{Higham2002}.

We consider the case where the matrix $\matrx{L}$ is stored in memory in a precision with unit roundoff $\varepsilon^{\store}$ and the triangular solves are computed in the same or higher precision with unit roundoff\footnote{The subscripts $\solve$ and $\store$ here stand for Solve and stoRe, respectively. They indicate that the corresponding $\varepsilon^{\solve}$- and $\varepsilon^{\store}$-precision are used for solving the triangular systems and for storing the matrix, respectively.} $\varepsilon^{\solve}$, $\varepsilon^{\solve} \leq \varepsilon^{\store}$. 
Storing $\matrx{L}$ in a lower precision than the precision used for solving the triangular systems may lead to faster memory accesses and thus to a faster runtime.
The smoothing routine is described in \Cref{alg:ic_smoothing}.

\begin{algorithm}[H]
\caption{IC smoother with zero initial approximation, $\matrx{ICS}(\vec{f})$.}
\label{alg:ic_smoothing}
\begin{algorithmic}[1]
\STATE{ $\vec{f} \leftarrow \mathrm{Round}(\vec{f},\varepsilon^{\solve}\text{-precision})$} 
\STATE{$\vec{v} \leftarrow \mathrm{Substitution}(\matrx{L},\vec{f})$}
\STATE{$\vec{w}\leftarrow\mathrm{Substitution}( \matrx{L}^{\top}, \vec{v} )$}
\RETURN{$\vec{w}$}
\end{algorithmic}
\end{algorithm}

Further, we present a finite precision error analysis of the application of the smoother. As mentioned in the introduction, we do not take into account the finite precision errors occurring when computing the factor $\matrx{L}$.
We first present a bound on the finite precision errors occurring when solving a general perturbed sparse triangular system via substitution.

\subsection{Finite precision error analysis of solving a perturbed sparse triangular system via substitution}

Let $\matrx{T}\in\R^{n \times n}$ be a sparse invertible triangular matrix and $\vec{b}\in \R^{n}$. We consider computing an approximate solution of 
$ \matrx{T}\vec{x} = \vec{b}$ 
using substitution. 
By modifying the proof of \cite[Theorem~8.5]{Higham2002}, using that there are a maximum of  $m_T$ nonzero entries in any row of $\matrx{T}$, we obtain the following result; we present its proof in \Cref{ape:sub_fp}.

\begin{lemma}\label{lemma:substitution}
Let $\matrx{T}$ and $\vec{b}$ be a matrix and a vector belonging to a finite precision arithmetic with unit roundoff $\varepsilon$ and $m_T\varepsilon<1$.
Let $\hat{\vec{x}}$ be the approximate solution of $\matrx{T}\vec{x}=\vec{b}$ computed via substitution in $\varepsilon$-precision. There exists a matrix $\matrx{E}$ such that
\begin{equation}\label{eq:backward_error_substitution}
    (\matrx{T} + \matrx{E})\hat{\vec{x}} = \vec{b},  \quad | \matrx{E} |\leq \varepsilon m_{T,\varepsilon}  | \matrx{T} |, \quad m_{T,\varepsilon} = m_T / (1-m_T\varepsilon).
\end{equation}
\end{lemma}

An important feature of this result is that the bound involves only the maximum number of non-zero elements in a row of the matrix $\matrx{T}$ but not the size of the matrix.
We use the result to derive a bound on the finite precision error when computing an approximate solution of $\matrx{T}\vec{x}=\vec{b}$ via substitution with the matrix $\matrx{T}$ rounded to a lower precision. 

\begin{lemma}\label{lemma:perturbed_substitution}
Let $\hat{\vec{x}}$ be the approximate solution of $\matrx{T} \vec{x} = \vec{b}$ computed via substitution in a finite precision arithmetic with unit roundoff $\varepsilon^{\solve}$, $m_T\varepsilon^{\solve}<1$, with the matrix $\matrx{T}$ rounded to a finite precision arithmetic with unit roundoff $\varepsilon^{\store}$, $\varepsilon^{\store}\geq\varepsilon^{\solve}$ and with $\vec{b}$ belonging to $\varepsilon^{\solve}$-precision arithmetic.
Let $\eta_{T}$ denote the constant $\varepsilon^{\store}  + \varepsilon^{\solve} m_{T,\varepsilon^{\solve}}   + \varepsilon^{\store}\varepsilon^{\solve} m_{T,\varepsilon^{\solve}}$.
If $\eta_{T} \underline{\kappa}_{T}<1 / 2 $, the Euclidean norm of the difference of $\vec{x}$ and $\hat{\vec{x}}$ can be bounded as 
\begin{align} \label{eq:perturbed_substitution_estimate_rhs_simlified}
\| \vec{x} - \hat{\vec{x}} \|
&\leq\eta_{T}  \underline{\kappa}_{T} (1 + 2\eta_{T} \underline{\kappa}_{T}) \| \vec{x}\|.
\end{align}
\end{lemma}
\begin{proof}
We first use the bound on the error when rounding a matrix to a lower precision and \Cref{lemma:substitution} to write a perturbed equation which the computed approximation $\hat{\vec{x}}$ satisfies. We then use this equation to derive the bound \eqref{eq:perturbed_substitution_estimate_rhs_simlified}.

Rounding the matrix $\matrx{T}$ to $\varepsilon^{\store}$-precision results in  $\matrx{T} + \Delta \matrx{T}$, where $ |\Delta \matrx{T}|\leq \varepsilon^{\store}| \matrx{T}|$; see the second inequality in \eqref{eq:fp_round}.
Using \Cref{lemma:substitution} for the perturbed problem $(\matrx{T} + \Delta \matrx{T} ) \tilde{\vec{x}} = \vec{b}$, there exists a matrix~$\matrx{F}$ such that 
\begin{equation}\label{eq:lemma:bound_F}
    \left( \matrx{T} + \Delta \matrx{T} + \matrx{F} \right) \hat{\vec{x}}  = \vec{b}, \quad
| \matrx{F} | \leq \varepsilon^{\solve} m_{T+\Delta T, \varepsilon^{\solve} }  | \matrx{T} + \Delta \matrx{T} |,
\quad m_{T+\Delta T,\varepsilon^{\solve}} = \frac{m_{T+\Delta T}}{1-m_{T + \Delta T}\varepsilon^{\solve}}.
\end{equation}  
Note that rounding a matrix can only result in it having fewer non-zero elements, and thus $m_{T+\Delta T}\leq m_T$ and consequently $m_{T+\Delta T,\varepsilon^{\solve}}\leq m_{T,\varepsilon^{\solve}}$.
From $\matrx{T}\vec{x}=\vec{b}$ and \eqref{eq:lemma:bound_F} we have
\begin{align}
\nonumber
\matrx{T} (\vec{x} - \hat{\vec{x}}) &= \vec{b} - \vec{b} + ( \Delta \matrx{T} + \matrx{F} )  \hat{\vec{x}} 
= ( \Delta \matrx{T} + \matrx{F} ) (\hat{\vec{x}}- \vec{x}) + ( \Delta \matrx{T} + \matrx{F} ) \vec{x},\\
\nonumber
\vec{x} - \hat{\vec{x}} &=   \matrx{T}^{-1} ( \Delta \matrx{T} + \matrx{F} ) (\hat{\vec{x}}- \vec{x}) + \matrx{T}^{-1} ( \Delta \matrx{T} + \matrx{F} ) \vec{x}, \\
\label{eq:lemma:x-hatx_2norm_bound}
\| \vec{x} - \hat{\vec{x}} \|  & \leq   \| \matrx{T}^{-1} \|   \| \Delta \matrx{T} + \matrx{F} \|  \|  \hat{\vec{x}} - \vec{x} \|  + \| \matrx{T}^{-1} \|  \| \Delta \matrx{T} + \matrx{F} \|  \| \vec{x} \|.
\end{align}
Using $| \Delta \matrx{T}|\leq \varepsilon^{\store} |\matrx{T}|$ and the bound in \eqref{eq:lemma:bound_F}, $| \Delta \matrx{T} + \matrx{F} |$ can be bounded as
\begin{align*}
| \Delta \matrx{T} + \matrx{F} | & \leq | \Delta  \matrx{T} | + | \matrx{F} | 
 \leq \varepsilon^{\store}| \matrx{T}| + \varepsilon^{\solve} m_{T,\varepsilon^{\solve}}  |\matrx{T} + \Delta \matrx{T} | 
 \leq \varepsilon^{\store}| \matrx{T}| + \varepsilon^{\solve}m_{T,\varepsilon^{\solve}} ( |\matrx{T} |  +  |\Delta \matrx{T} |  )\\
 & \leq (\varepsilon^{\store}  + \varepsilon^{\solve} m_{T,\varepsilon^{\solve}}  
+ \varepsilon^{\store}\varepsilon^{\solve} m_{T,\varepsilon^{\solve}}) |\matrx{T} | = \eta_{T} |\matrx{T} |.
\end{align*}
This yields (see, e.g., \cite[Lemma~6.6,~case~(b)]{Higham2002}) the bound 
$\| \Delta \matrx{T} + \matrx{F} \| \leq \eta_{T} \| |\matrx{T} | \|$.
Using this in \eqref{eq:lemma:x-hatx_2norm_bound} and using the definition of $\underline{\kappa}_{T}$ leads to
\begin{equation*}
\| \vec{x} - \hat{\vec{x}} \| 
 \leq  \| \matrx{T}^{-1} \|  \eta_{T} \|  | \matrx{T} |  \|  \|  \hat{\vec{x}} - \vec{x} \|  + \| \matrx{T}^{-1} \| \eta_{T}  \|  | \matrx{T} |  \|   \| \vec{x} \| 
 =   \eta_{T} \underline{\kappa}_{T}   \|  \hat{\vec{x}} - \vec{x} \|  + \eta_{T}  \underline{\kappa}_{T} \| \vec{x} \|,
\end{equation*}
and subsequently, 
$(1 - \eta_{T}  \underline{\kappa}_{T})\| \vec{x} - \hat{\vec{x}} \|
\leq  \eta_{T}  \underline{\kappa}_{T} \| \vec{x}\|$. 
Using the assumption $ \eta_{T}  \underline{\kappa}_{T}< 1 / 2 < 1 $ and dividing by $1 - \eta_{T}  \underline{\kappa}_{T}$ gives
\begin{equation*}
\| \vec{x} - \hat{\vec{x}} \|
\leq  \eta_{T}  \underline{\kappa}_{T} 
(1 - \eta_{T} \underline{\kappa}_{T})^{-1}
\| \vec{x}\|.
\end{equation*}
Using that $(1-t)^{-1} \leq 1+2t$ for $t \in (0,1/2 )$ 
and the assumption $ \eta_{T}   \underline{\kappa}_{T}<1 / 2$ we can bound 
$(1-\eta_{T} \underline{\kappa}_{T})^{-1}$ by $1+2  \eta_{T} \underline{\kappa}_{T}$ and obtain \eqref{eq:perturbed_substitution_estimate_rhs_simlified}. 
\end{proof}

\subsection{Finite precision error analysis of mixed precision IC smoother}

In this section we use the results from the previous section and present a bound on the finite precision error in the application of the IC smoother.

\begin{theorem}\label{thm:IC_fp_error}
Let $\vec{w}$ and $\hat{\vec{w}}$ be the approximations computed by applying the IC smoother (\Cref{alg:ic_smoothing}) to a vector $\vec{f}$ in exact arithmetic (without rounding the matrix $\matrx{L}$ to a lower precision) and in finite precision, respectively.
Let $\underline{m}_{L,\varepsilon^{\solve}} = \underline{m}_{L} / (1-\underline{m}_{L} \varepsilon^{\solve})$
and $\eta_{L} = \varepsilon^{\store}  + \varepsilon^{\solve}\underline{m}_{L,\varepsilon^{\solve}}   + \varepsilon^{\store}\varepsilon^{\solve}\underline{m}_{L,\varepsilon^{\solve}}$.
Assuming $\underline{m}_{L} \varepsilon^{\solve}<1$ and $  \eta_{L} \underline{\kappa}_{L}<1 / 2$ the Euclidean norm of the finite precision error $\delta_{\IC}=\vec{w}-\hat{\vec{w}}$ can be bounded as
\begin{equation}\label{eq:thm:rounding_error_smoother}
\| \delta_{\IC} \| \leq \Lambda_{\IC} \| \vec{f} \|, \quad \Lambda_{\IC} = (2(\varepsilon^{\store}  + \varepsilon^{\solve}(\underline{m}_{L,\varepsilon^{\solve}}+1 / 2 )) \underline{\kappa}_{L} + R)
\|  \matrx{L}^{-1} \|^2 ,
\end{equation}
where the remainder $R$ contains higher order terms, i.e., terms which involve $(\varepsilon^{\solve})^2$, $(\varepsilon^{\store})^2$, or $\varepsilon^{\solve}\varepsilon^{\store}$.
\end{theorem}

\begin{proof}
The terms $\vec{v}$ and $\vec{w}$ computed in exact arithmetic satisfy $\vec{v}=\matrx{L}^{-1}\vec{f}$ and $\vec{w}=\matrx{L}^{-\top}\matrx{L}^{-1}\vec{f}$.
Let $\hat{\vec{f}}$, $\hat{\vec{v}}$, and $\hat{\vec{w}}$ denote the corresponding terms computed in finite precision. 
Below, $R$, $R_1$, and $R_2$ denote remainders which contain only higher order terms, i.e., terms which involve $(\varepsilon^{\solve})^2$, $(\varepsilon^{\store})^2$, or $\varepsilon^{\solve}\varepsilon^{\store}$.

Rounding $\vec{f}$ to $\varepsilon^{\solve}$-precision arithmetic results in $\hat{\vec{f}}=\vec{f}+\delta_f$, where using \eqref{eq:fp_round}, $\| \delta_f \| \leq \varepsilon^{\solve}\| \vec{f} \|$.
Computing the first substitution in finite precision results in $\hat{\vec{v}} = \matrx{L}^{-1}  (\vec{f} + \delta_f) + \delta_{v} $, where, using \Cref{lemma:perturbed_substitution} and inequality \eqref{eq:perturbed_substitution_estimate_rhs_simlified},
\begin{align*}
\| \delta_v \| 
& \leq
 \eta_{L}   \underline{\kappa}_{L}
(1+2\eta_{L} \underline{\kappa}_{L} )\|   \matrx{L}^{-1}  (\vec{f} + \delta_f)  \| 
\leq 
 \eta_{L}   \underline{\kappa}_{L}
(1+2\eta_{L} \underline{\kappa}_{L} )(1+\varepsilon^{\solve}) \|   \matrx{L}^{-1} \| \|  \vec{f} \| 
\\
&\leq 
\underbrace{((\varepsilon^{\store}  + \varepsilon^{\solve}\underline{m}_{L,\varepsilon^{\solve}}) \underline{\kappa}_{L} + R_1)}_{K_1}
\| \matrx{L}^{-1} \| \|  \vec{f} \|.
\end{align*}
The term  $\hat{\vec{v}}$ can be written as $\hat{\vec{v}} = \vec{v} + \Delta_v$, where $\Delta_v = \matrx{L}^{-1}  \delta_f + \delta_v$ and
$\| \Delta_v\| \leq (K_1+\varepsilon^{\solve}) \| \matrx{L}^{-1} \| \|  \vec{f} \|$.

Computing the second substitution in finite precision results in 
$\hat{\vec{w}}=\matrx{L}^{-\top}(\vec{v} + \Delta_v ) + \delta_{w}$, where using \Cref{lemma:perturbed_substitution} and inequality \eqref{eq:perturbed_substitution_estimate_rhs_simlified},
\begin{align*}
\| \delta_w \|
& \leq   \eta_{L}   \underline{\kappa}_{L}(1+2 \eta_{L} \underline{\kappa}_{L})
\| \matrx{L}^{-\top}(\vec{v} + \Delta_v )\|
 \\
&\leq 
\eta_{L}  \underline{\kappa}_{L} 
(1 +2 \eta_{L} \underline{\kappa}_{T})
 (1+ K_1 + \varepsilon^{\solve}) \|  \matrx{L}^{-1}  \|^2 \|  \vec{f} \|  \\
&\leq 
\underbrace{((\varepsilon^{\store}  + \varepsilon^{\solve}\underline{m}_{L,\varepsilon^{\solve}}) \underline{\kappa}_{L} + R_2)}_{K_2}
\| \matrx{L}^{-1} \|^2 \| \vec{f} \|.
\end{align*}
The term $\hat{\vec{w}}$ can then be written as $\hat{\vec{w}} = \vec{w} + \Delta_w$, where $\Delta_w = \delta_w + \matrx{L}^{-1}\Delta_v$ is the accumulated error and 
\begin{equation*}
\| \Delta_w \|
\leq ( K_2  +  K_1 + \varepsilon^{\solve}) \| \matrx{L}^{-1} \|^{2} \| \vec{f} \| \leq 
(2(\varepsilon^{\store}  + \varepsilon^{\solve}\underline{m}_{T,\varepsilon^{\solve}}) \underline{\kappa}_{L} + \varepsilon^{\solve} + R) \|  \matrx{L}^{-1}  \|^2 \|  \vec{f} \|.
\end{equation*}
Since $\underline{\kappa}_T\geq1$, the term $\varepsilon^{\solve}$ can be bounded by $\varepsilon^{\solve}\underline{\kappa}_T$ which yields the bound \eqref{eq:thm:rounding_error_smoother}.
\end{proof}

We remark that the bound \eqref{eq:thm:rounding_error_smoother} is the worst case scenario bound and the actual error could be significantly smaller.
The number $\underline{m}_L$ depends on the sparsity pattern of the matrix $\A$ and the fill-in that occurs in the IC factorization. We see that the requirements on the $\varepsilon^{\store}$- and  $\varepsilon^{\solve}$-precisions differ in the multiplicative constant $\underline{m}_{L,\varepsilon^{\solve}}+ 1 / 2$.

\section{Mixed precision V-cycle with IC smoothing}\label{sec:IC-V-cycle}
We summarize the results on the effects of finite precision errors in the V-cycle (\Cref{alg:V-cycle}) with IC smoothing (\Cref{alg:ic_smoothing}) in this section.

We consider the case where IC smoothing is used on all fine levels of the V-cycle.
In the notation of \Cref{sec:V-cycle} we have $\matrx{M}_j =\matrx{L}^{-\top}_j  \matrx{L}^{-1}_j $, $j=1,\ldots,J$, where $ \matrx{L}_j \matrx{L}^{\top}_j$ is an IC factorization of $\A_j$.
The precisions used on the $j$th level are:  
\vspace{1mm}
\begin{itemize}
\item[] $\dot{\varepsilon}_j$-precision, for computing residual, restriction, prolongation, and correction,
\item[] $\varepsilon^{\store}_j$-precision, for storing the matrix $\matrx{L}_j$ in the memory,  and
\item[] $\varepsilon^{\solve}_j$-precision, for solving triangular systems with $\matrx{L}_j$ and $\matrx{L}^{\top}_j$ via substitution.
\end{itemize}
\vspace{1mm}
We assume that the precisions used for the smoothing are lower or equal to the precision used for the other operations, i.e., $\varepsilon^{\store}_j \geq \varepsilon^{\solve}_j \geq \dot{\varepsilon}_j$.
 
Assuming that the $\varepsilon^{\store}_j$- and $\varepsilon^{\solve}_j$- precisions are chosen such that the assumptions of  \Cref{thm:IC_fp_error} are satisfied, the theorem shows that assumption \eqref{eq:smoother_assumption_fp_V-cycle}  on the finite precision error when applying the smoother on the $j$th level is satisfied, giving
\begin{equation}
\label{eq:lambda_M_estimate}
\matrx{M}_j \vec{f}_j + \delta_{M_j},\quad \| \delta_{M_j} \| \lesssim \Lambda_{M_j}  \| \vec{f}_j\|, \quad \Lambda_{M_j}  = 
 2(\varepsilon^{\store}_j + \varepsilon^{\solve}_j(\underline{m}_{L_j,\varepsilon^{\solve}_j} + 1 / 2)
) \underline{\kappa}_{L_j} 
\|  \matrx{L}^{-1}_j  \|^2;
\end{equation}
we use $\lesssim$ to indicates that we have dropped the remainder containing higher order terms.
We note that we are not able to theoretically verify the assumption \eqref{eq:V_cycle_smoothing_A_norm_assumption}. It can, however, be verified  numerically in concrete settings.

If the other assumptions of \Cref{thm:V-cycle} are satisfied, we get the following bound on the relative finite precision error after the V-cycle (\Cref{alg:V-cycle}) application:
\begin{align}\label{eq:V-cycle-IC-fp}
\begin{split}
\frac{\|  \vec{y}_{\V} - \hat{\vec{y}}_{\V} \|_{\A}}{\| \vec{y} \|_{\A}}  \lesssim 
\Lambda_0  & + 3 \sum^{J}_{j=1}   (\varepsilon^{\store}_j + \varepsilon^{\solve}_j(\underline{m}_{L_j,\varepsilon^{\solve}_j} + 1 / 2)
) \underline{\kappa}_{L_j}  
\|  \matrx{L}^{-1}_j  \|^2 \| \A_j \| 
\\ &  
\quad +  \sum^{J}_{j=1} \dot{\varepsilon}_j \kappa^{\frac{1}{2}}_{A_j}  (C_{1,j}   \| \matrx{L}^{-1}_j \|^2 + C_{2,j} ),
\end{split}
\end{align}
where the constants $C_{1,J}$ and $C_{2,J}$ depends only on  $\| \matrx{A}_j\|$, $ \| | \matrx{A}_j | \|$, $\| \matrx{P}_j\|$, $ \| | \matrx{P}_j | \|$, $m_{A_j,\dot{\varepsilon}_j}$,  $\underline{m}_{P_j,\dot{\varepsilon}_j}$, and the ratio $\| \A^{-1}_{j-1}\|^{\frac{1}{2}} / \| \A^{-1}_j \|^{\frac{1}{2}}$. 
If the $\varepsilon_j$-, $\varepsilon^{\store}_j$-, and $\varepsilon^{\solve}_j$-precisions are chosen such that the right-hand side of \eqref{eq:V-cycle-IC-fp} is much smaller than one, the relative finite precision error is much smaller than one and the contraction factor of the V-cycle should not be significantly affected by finite precision errors. We see that the upper bound on the finite precision error of the coarsest-level solver $\Lambda_0$ is present as a standalone term.
The requirements given by the estimate \eqref{eq:V-cycle-IC-fp} on the $\varepsilon^{\store}_j$- and $\varepsilon^{\solve}_j$-precisions depend on properties of the IC factor $\matrx{L}_j$ ($\underline{\kappa}_{L_j}$, $\|  \matrx{L}^{-1}_j  \|^2$ and $\underline{m}_{L_j}$) and on the norm of $\A_j$,   $\| \A_j \|$,  but not, for example, on the square root of its condition number $\kappa^{\frac{1}{2}}_{A_j}$.
The requirements on the $\dot{\varepsilon}_j$-precision depend on the square root of the condition number of $\A_j$, $\kappa^{\frac{1}{2}}_{A_j}$, the property of the IC factor $\matrx{L}_j$, $\|  \matrx{L}^{-1}_j  \|^2$, and the constants $C_{1,j}$ and $C_{2,j}$.
This allows us to make the following observations:
\begin{itemize}
\item The $\dot{\varepsilon}_j$-, $\varepsilon^{\store}_j$-, and $\varepsilon^{\solve}_j$-precisions can be chosen differently on different levels, based on the properties of the multigrid hierarchy. 

\item If $\kappa^{\frac{1}{2}}_{A_j}$ (times the corresponding terms) is decreasing when going from fine to coarse levels, which is typical in multgrid settings (see, e.g., \cite{Tamstorf2021} or the experiments presented below), the $\dot{\varepsilon}_j$-precision can be chosen to be progressively lower when going from fine to coarse levels. This is in line with the results in \cite{McCormick2021,Tamstorf2021,McCormick2023}.
\item In settings where $\underline{\kappa}_{L_j}$ (times the corresponding terms) is smaller than $\kappa^{\frac{1}{2}}_{A_j}$ (times the corresponding terms), the $\varepsilon^{\store}_j$- and $\varepsilon^{\solve}_j$-precisions can be chosen to be lower than the $\dot{\varepsilon}_j$-precision.  
\item The requirements on the $\varepsilon^{\store}_j$-precision and $\varepsilon^{\solve}_j$-precision differ in the multiplicative factor containing $\underline{m}_{L_j}$. If this term is sufficiently large the $\varepsilon^{\store}_j$-precision can be chosen lower than the $\varepsilon^{\solve}_j$-precision.

\item If $\underline{\kappa}_{L_j}$ (times the corresponding terms) is decreasing when going from fine to coarse levels the $\varepsilon^{\store}_j$- and $\varepsilon^{\solve}_j$-precisions can be chosen to be progressively lower when going from fine to coarse levels.
In settings where $\underline{\kappa}_{L_j}$ (times the corresponding terms) are nearly the same on different levels, the $\varepsilon^{\store}_j$- and $\varepsilon^{\solve}_j$-precisions can be chosen to be the same on these levels.

\item Using IC smoothing with different fill-ins yields different factors $\matrx{L}_j$ with different properties and thus different requirements not only on the smoothing-related  $\varepsilon^{\store}_j$- and $\varepsilon^{\solve}_j$-precisions, but also on the $\dot{\varepsilon}_j$-precisions. 
\item The precision used for applying the coarsest-level solver should be chosen based on the concrete coarsest-level problem and solver such that the resulting coarsest-level finite precision error is small. Certain settings may require the precision used for applying the coarsest-level solver to be higher than the $\dot{\varepsilon}_j$-precisions on some finer level.
\end{itemize}
We illustrate these findings in numerical experiments in \Cref{sec:num-exp-simul,sec:num-exp-gingko}.

\section{Scaling system matrices and right-hand sides}\label{sec:scaling}
Rounding matrices or vectors to low precision and computing in low precision can result in overflow or underflow errors; we refer, e.g., to the discussions in \cite{Higham2019}. Scaling the data before rounding and computing can help to partially overcome this issue.

In the experiments below we use a simple scaling for the system and prolongation matrices in a multigrid hierarchy which preserves the Galerkin condition.
The matrix $\A_j$ on the $j$th level, $j=0,\ldots,J$, is scaled as $\bar{\A}_j = s_j \A_j$, where $s_j = 1 / \max_{k,\ell} |[\A_j]_{k,\ell}|$. The prolongation matrix on the $j$th level, $j=0,\ldots,J-1$, is scaled as $\bar{\matrx{P}}_j = \frac{\sqrt{s_{j-1}}}{\sqrt{s_{j}}}\matrx{P}_j$.
The Galerkin condition then holds (in exact arithmetic) also for the scaled matrices
\begin{equation*}
\bar{\matrx{P}}_j^\top \bar{\A}_{j} \bar{\matrx{P}}_j =
\frac{\sqrt{s_{j-1}}}{\sqrt{s_{j}}} \matrx{P}_j^\top s_j\A_{j} 
\frac{\sqrt{s_{j-1}}}{\sqrt{s_{j}}} \matrx{P}_j
=s_{j-1} \matrx{P}_j^\top \A_{j} \matrx{P}_j
= s_{j-1} \A_{j-1} = \bar{\A}_{j-1}.
\end{equation*}

Scaling can be also applied to a right-hand side vector before calling a smoothing routine; see e.g., \cite[Section~6]{Carson2018}. 
We first compute the infinity norm of the right-hand side vector $\vec{f}$, i.e.,  $s_f=\| \vec{f}\|_{\infty}$. The right-hand side is then scaled as $\bar{\vec{f}} = s^{-1}_f \vec{f}$ and the smoothing is called with the scaled vector $\bar{\vec{f}}$. The result is afterwards re-scaled back by multiplying with $s_f$.

We remark that the discussed scaling may help with staying inside the range of a low precision arithmetic; however, it is not guaranteed.

\section{Numerical experiment with simulated floating point arithmetics}
\label{sec:num-exp-simul}
The goal of this experiment is to illustrate the theoretical findings summarized at the end of \Cref{sec:IC-V-cycle}.
We solve systems of linear equations obtained by discretization of the 1D elliptic PDE:
 find $u:(0,1) \rightarrow \R$ such that 
\begin{equation*}
    - u'' = f \quad \textrm{in } (0,1), \quad u(0)  = u(1) =  0,
\end{equation*}
where $f$ is chosen to correspond to the solution $u(x) = x(x-1) \sin(2\pi x)$.

The problem is discretized using the continuous Galerkin FE method with piecewise polynomials of degree five (FEM-P5) on a hierarchy of $15$ uniformly refined meshes.
We consider this 1D problem since it allows us to work with multigrid hierarchies with a large number of levels. This benefits the illustration of how the requirements on the precisions change on different levels. 
The matrices are assembled in the finite element software FEniCS \cite{Alnaes2015} in double precision. 
We modify the system matrices, the prolongation matrices, and the right-hand side vectors so that the resulting systems contain just free-node variables. The Galerkin condition is then satisfied on all coarse levels.
We scale the system and prolongation matrices and the right-hand side vectors using the strategy described in \Cref{sec:scaling}.
We also filter values of the system matrices and the prolongation matrices after scaling at the level $5\cdot 10^{-16}$ and $5\cdot 10^{-12}$, respectively. 
The numbers of degrees of freedom (DoF) grows approximately by a factor of two with each fine level from $24$ to $409,599$. The data and codes for reproducing the experiments in this paper can be found
at \url{https://doi.org/10.5281/zenodo.13858606}.

We consider the case where the coarsest level is fixed and we solve the problems $\A_J \vec{x}_J=\vec{b}_J$, $J=2,\ldots,14$, using the IR method where the inner solver is the geometric V-cycle method (\Cref{alg:V-cycle}) with IC smoothing (\Cref{alg:ic_smoothing}) and $J+1$ levels (referred to as IR-V-cycle-IC).
The computation is done in MATLAB.
The residual computation and the approximate solution update within IR are done in double precision.
The method is run starting with zero initial approximation and stopped when the absolute algebraic error in the $\A_J$-norm is (approximately) less than $10^{-5}$ (the reference solution for computing the algebraic error is approximated by the MATLAB backslash operator in double precision).
The initial algebraic error is approximately $10^{-1}$. 

We note that in this experiment, we focus only on the effects of the finite precision errors on the algebraic error. It is possible that in order to compute an accurate solution of the PDE, the matrices on fine levels would have to be assembled in a precision higher than double and the residual computation and solution update in IR would also have to be done in a higher precision; see the experiments in \cite{Tamstorf2021}. Being aware of the limitations of our setting, we believe that the experiment well illustrates the theoretical results on the effects of finite precision errors on the algebraic error.

We use V-cycles where the $\dot{\varepsilon}_j$-precisions
are the same on all fine levels, i.e., $\dot{\varepsilon}_j=\dot{\varepsilon}_J$, $j=1,\ldots,J$. 
The same holds for the precisions used in IC smoothers, where we additionally assume that the $\varepsilon^{\store}_j$-precision is the same as the $\varepsilon^{\solve}_j$-precision, i.e.,  $\varepsilon^{\store}_j=\varepsilon^{\solve}_j=\varepsilon^{\solve}_J$, $j=1,\ldots,J$.
We consider two variants of IC smoothers: IC($0$), based on factorization with zero fill-in, and ICT(dpt=$5 \cdot 10^{-3}$), with a local dropping tolerance $5 \cdot 10^{-3}$.
Allowing fill-in in the IC factorization typically yields a better approximation of the matrix $\A_j$ by $\matrx{L}_j\matrx{L}^{\top}_j$ and consequently to a more effective smoothing routine.
The solver on the coarsest-level, the MATLAB backslash operator, is applied in double precision to a problem with matrix $\A_0$ rounded to $\dot{\varepsilon}_J$-precision, and the computed coarsest-level approximation is subsequently rounded to $\dot{\varepsilon}_J$-precision.

We use simulated floating point arithmetic in the experiment
since it enables us to compute in multiple arithmetics of mildly varying precisions, which benefits the illustrations. In particular, we utilize the Advanpix toolbox \cite{advanpix}. This toolbox requires specifying the number of decimal digits $d$, simulating\footnote{The toolbox has $64$ bits for representing the exponent, except for the variant with $d=34$ where it is $15$ bits; see e.g., \cite[Section~8]{Tamstorf2021}. The large number of bits for representing the exponent results in the computation not being affected by the limited range as it is when using the standard single and especially half precision, which have $8$ and $5$ bits for storing the exponent, respectively.} the floating point arithmetic with approximate unit roundoff $10^{-d}$.

We first run the computation with all the precisions set to double precision. 
Then we assume that the  $\dot{\varepsilon}_J$- and $\varepsilon^{\solve}_J=\varepsilon^{\store}_J$-precisions are the same and we run the computation using the Advanpix toolbox simulating $\dot{\varepsilon}_J=\varepsilon^{\solve}_J=\varepsilon^{\store}_J=10^{-\dot{d}_J}$, for $\dot{d}_J=1,2,\ldots$.
We find the smallest $\dot{d}_J$, denoted as $\dot{d}_{J,\mathrm{min}}$, for which the method converges in the same number of IR iterations as the corresponding variant in double precision.
Further, we fix the $\dot{\varepsilon}_J$-precision as $\dot{\varepsilon}_J=10^{\dot{d}_{J,\mathrm{min}}}$ and run the experiments simulating $\varepsilon^{\solve}_J=\varepsilon^{\store}_J=10^{-d^{\solve}_J}$, $d^{\solve}_J = 1,2,\ldots$ . We again find the minimal $d^{\solve}_J$, denoted as $d^{\solve}_{J,\mathrm{min}}$, for which the method converges in the same number of IR iterations as the corresponding variant in double precision.

\subsubsection*{Expectations based on the finite precision error analysis}
From the finite precision error analysis of IR in the $\A$-norm \cite[Sections~3-4]{McCormick2021}, we know that the convergence of IR is determined by the precisions chosen for computing the residual and the approximate solution update and by the contraction factor of the inner solver applied in finite precision.  
Since we are comparing variants of IR 
that only differ in terms of the precision formats used for the inner solver, its enough to discuss how these precisions affect the contraction factor of the inner solver---here one V-cycle starting with zero initial approximation.
We will use the results presented in \Cref{sec:IC-V-cycle}.

We approximately compute the values of the terms in the estimate \eqref{eq:V-cycle-IC-fp} of the relative finite precision error after applying the V-cycle (\Cref{alg:V-cycle}). The values of the square root of the condition number of the system matrices $\kappa^{\frac{1}{2}}_{A_j}$, the condition number of the IC factors  $\underline{\kappa}_{L_j}$, and the norms  $\|\matrx{L}_j^{-1} \|^2$ are plotted in \Cref{fig:mat_properties_1D_and_sizes1D_3D} (left).

We see that $\kappa^{\frac{1}{2}}_{A_{j}}$ grows approximately by a factor of two with each finer level, i.e., $2\kappa^{\frac{1}{2}}_{A_j} \approx \kappa^{\frac{1}{2}}_{A_{j+1}} $.
The values of $\underline{\kappa}_{L_j}$ and $\|\matrx{L}_j^{-1} \|^2$ do not significantly change on different  levels. The values of $\|\matrx{L}^{-1}_j \|^2$ and $\underline{\kappa}_{L_j}$ are larger for the variant with ICT(dpt=$5 \cdot 10^{-3}$) than the corresponding values for the variant with IC(0).
We also approximately compute the following properties (they are nearly the same on all levels) $\| \A_j \|\approx \| |\A_j | \|\approx2.6$,~ $ \| \A^{-1}_{j-1}\|^{\frac{1}{2}} / \| \A^{-1}_j \|^{\frac{1}{2}}\approx 2$,
$\| \matrx{P}_j \| \approx 3.2$,~ $ \| |\matrx{P}_j| \|\approx3.6$, $m_{A_j}=11$, $\max_j\underline{m}_{P_j}=12$, $\underline{m}_{L_j}=10$ for both IC(0) and ICT(dpt=$5 \cdot 10^{-3}$); we note that $\underline{m}_{L_j}$ is the maximum number of nonzero entries in either a row or a column of $\matrx{L}_j$.

\begin{figure}
\centering
\setlength\figureheight{4.5cm}
\setlength\figurewidth{0.99\textwidth}
\definecolor{mycolor1}{RGB}{94, 201, 98}%
\definecolor{mycolor2}{RGB}{33, 145, 140}%
\definecolor{mred}{RGB}{237, 121, 83}%
\definecolor{mblue}{RGB}{65, 68, 135}%
\begin{tikzpicture}
\begin{axis}[%
width=0.35\figurewidth,
height=\figureheight,
at={(0.25\figurewidth,0\figureheight)},
xlabel={level $j$},
ylabel = {},
ylabel style={yshift=-1.2em},
xmin = 1,
xmax = 14,
ymode = log,
xtick = {1,5,10,14},
xticklabels = {1,5,10,14},
yminorticks=true,
legend pos= outer north east,
legend style = {draw = none},
legend cell align= left,
axis background/.style={fill=white}
]
\addplot [color=mycolor1, style = solid, line width=1.0pt, mark=asterisk, mark options={solid}]
  table[row sep=crcr]{%
1	3.75E+01	\\
2	7.50E+01	\\
3	1.50E+02	\\
4	3.00E+02	\\
5	6.01E+02	\\
6	1.20E+03	\\
7   2.40E+03    \\
8	4.81E+03	\\
9	9.61E+03	\\
10	1.92E+04	\\
11	3.84E+04	\\
12  7.69E+04    \\
13	1.54E+05	\\
14	3.08E+05	\\
15  6.15E+05    \\
};
\label{line:sqrt_cond_A_1D}

\addplot [color=mred, style = dashed, line width=1.0pt, mark=+, mark options={solid}]
  table[row sep=crcr]{%
1	1.08E+01	\\
2	1.08E+01	\\
3	1.08E+01	\\
4	1.08E+01	\\
5	1.08E+01	\\
6   1.08E+01    \\
7	1.08E+01	\\
8	1.08E+01	\\
9	1.08E+01	\\
10	1.08E+01	\\
11  1.08E+01    \\
12	1.08E+01	\\
13	1.08E+01	\\
14  1.08E+01    \\
};
\label{line:cond_L_ic(0)_1D}

\addplot [color=mred, style = solid, line width=1.0pt, mark=o, mark options={solid}]
  table[row sep=crcr]{%
1	3.98E+01	\\
2	3.99E+01	\\
3	3.99E+01	\\
4	3.99E+01	\\
5	3.99E+01	\\
6   3.99E+01    \\
7	3.99E+01	\\
8	3.99E+01	\\
9	3.99E+01	\\
10	3.99E+01	\\
11  3.99E+01    \\
12	3.99E+01	\\
13	3.99E+01	\\
14  3.99E+01    \\
};
\label{line:norm_L_sq_ic(0)_1D}

\addplot [color=mblue, style = dashed, line width=1.0pt, mark=x, mark options={solid}]
  table[row sep=crcr]{%
1	4.34E+01	\\
2	4.67E+01	\\
3	4.76E+01	\\
4	4.78E+01	\\
5	4.79E+01	\\
6   4.79E+01    \\
7	4.79E+01	\\
8	4.79E+01	\\
9	4.79E+01	\\
10	4.79E+01	\\
11  4.79E+01    \\
12	4.79E+01	\\
13	4.79E+01	\\
14  4.79E+01    \\
};
\label{line:cond_L_ic(dpt)_1D}

\addplot [color=mblue, style = solid, line width=1.0pt, mark=square, mark options={solid}]
  table[row sep=crcr]{%
1	4.80E+02	\\
2	5.50E+02	\\
3	5.70E+02	\\
4	5.76E+02	\\
5	5.77E+02	\\
6   5.77E+02    \\
7	5.77E+02	\\
8	5.77E+02	\\
9	5.77E+02	\\
10	5.77E+02	\\
11  5.77E+02    \\
12	5.77E+02	\\
13	5.77E+02	\\
14  5.77E+02    \\
};
\label{line:norm_L_sq_ic(dpt)_1D}
\legend{ $\kappa^{\frac{1}{2}}_{A_j}$, {$\underline{\kappa}_{L_j}$ (IC(0))}, $\underline{\kappa}_{L_j}$ (ICT), $\| \matrx{L}^{-1}_j \|^2$ (IC(0)), $\| \matrx{L}^{-1}_j \|^2$ (ICT)}
\end{axis}
\end{tikzpicture}
\caption{1D Poisson problem, FEM-P5 disc. Properties of $\A_j$ and $\matrx{L}_j$ for IC(0) and  \mbox{ICT(dpt=$5 \cdot 10^{-3} $)}.
} 
\label{fig:mat_properties_1D_and_sizes1D_3D}
\end{figure}

Since the coarsest-level solver is applied in double precision to a small well-conditioned problem, we expect the associated finite precision error to be negligible.
We note that we assume that the precisions are fixed on all levels, i.e., $\varepsilon^{\store}_j=\varepsilon^{\solve}_j=\varepsilon^{\solve}_J$ and $\dot{\varepsilon}_j=\dot{\varepsilon}_J$. 
The estimate \eqref{eq:V-cycle-IC-fp} can be simplified as
\begin{align*}
\begin{split}
\frac{\|  \vec{y}_{\V,J} - \hat{\vec{y}}_{\V,J} \|_{\A_J}}{\| \vec{y}_J \|_{\A_J}} & \lesssim 6 \| \A_J \|  
( \underline{m}_{L_J,\varepsilon^{\solve}_J} + 3 / 2))
\| \matrx{L}^{-1}_J \|^2 
\varepsilon^{\solve}_J\underline{\kappa}_{L_J} J \\ & 
\quad
\ +   2(C_{1,J}  \| \matrx{L}^{-1}_J \|^2 + C_{2,J} )  \dot{\varepsilon}_J  \kappa^{\frac{1}{2}}_{A_J}.
\end{split}
\end{align*}
We also used that the sum  $\sum^{J}_{j=1} \kappa^{\frac{1}{2}}_{A_j}$ can be approximated by $2 \kappa^{\frac{1}{2}}_{A_J}$, which is a consequence of  $\kappa^{\frac{1}{2}}_{A_j} \approx \kappa^{\frac{1}{2}}_{A_{j+1}} / 2$. 

Since $\underline{\kappa}_{L_J}J$ grows linearly with an increasing number of levels $J$ ($\underline{\kappa}_{L_J}$ is nearly constant), whereas $\kappa^{\frac{1}{2}}_{A_J}$ grows approximately exponentially as $\kappa^{\frac{1}{2}}_{A_0}2^J$ with increasing number of levels $J$, we expect that for a larger number of levels $J=8,\ldots,14$, the $\varepsilon^{\solve}_J$-precision could be significantly lower than the $\dot{\varepsilon}_J$-precision while preserving the same convergence rate of the IR-V-cycle method. This is valid for both the IC($0$) and ICT(dpt=$5 \cdot 10^{-3}$) variants.

Since the values of $\|\matrx{L}^{-1}_J \|^2$ and $\underline{\kappa}_{L_J}$ are larger for the variant with ICT(dpt=$5 \cdot 10^{-3}$) than the corresponding values for the variant IC(0), we expect that the $\dot{\varepsilon}_J$- and $\varepsilon^{\solve}_J$-precisions for the variant with ICT(dpt=$5 \cdot 10^{-3}$) may have to be chosen higher than for the IC(0) variant. 

\subsubsection*{Results}
The computed values of $\dot{d}_{J,\mathrm{min}}$ and $d^{\solve}_{J,\mathrm{min}}$ are summarized in \Cref{fig:1D_IC_results} together with the number of IR iterations required to reach the chosen accuracy.
We see that the variant with ICT(dpt=$5\cdot 10^{-3}$) requires significantly fewer IR iterations than the variants with IC(0). Regardless of the variant of the IC factorization, the values of $d^{\solve}_{J,min}$ corresponding to the $\varepsilon^{\solve}_J=\varepsilon^{\store}_J$-precision used in the smoothing are smaller than the corresponding values of $\dot{d}_{J,min}$ corresponding to the $\dot{\varepsilon}_J$-precision used for computing the residual, prolongation, restriction, and correction inside the V-cycle. Moreover $\dot{d}_{J,min}$ increases when increasing $J$, while $d^{\solve}_{J,min}$ stays constant. This illustrates that the $\varepsilon^{\solve}_J=\varepsilon^{\store}_J$-precision may be, in some settings, significantly lower than the $\dot{\varepsilon}_J$-precision.
We observe that the values of $\dot{d}_{J,\mathrm{min}}$  and $d^{\solve}_{J,\mathrm{min}}$ for the variant with ICT(dpt=$5\cdot10^{-3}$) are larger than or equal to the corresponding values for the variant with IC(0). In other words the variant with ICT(dpt=$5\cdot10^{-3}$) requires higher or equal $\dot{\varepsilon}_J$- and $\varepsilon^{\solve}_J=\varepsilon^{\store}_J$-precisions than the variant with IC(0).
We conclude that the results align with the theoretical derivations.

Even though we run the experiments with the $\dot{\varepsilon}_J$-precision fixed for all levels, this experiment also illustrates that the $\dot{\varepsilon}_j$-precision, $j=1\ldots,J$, could be chosen to be lower on the coarse levels and progressively increased.

\begin{figure}
\centering
\setlength\figureheight{6cm}
\setlength\figurewidth{0.99\textwidth}
\definecolor{mycolor1}{RGB}{94, 201, 98}%
\definecolor{mycolor2}{RGB}{33, 145, 140}%
\begin{tikzpicture}

\begin{axis}[%
width=0.49\figurewidth,
height=\figureheight,
at={(0\figurewidth,0\figureheight)},
xlabel={level $J$},
ylabel = { number of  digits},
ylabel style={yshift=-1.2em},
xmin = 2,
xmax = 14,
ytick = {2,3,4,5,6,7,16},
yticklabels = {2,3,4,5,6,7,16},
xtick = {2,3,4,5,6,7,8,9,10,11,12,13,14},
xticklabels = {2,3,4,5,6,7,8,9,10,11,12,13,14},
yminorticks=true,
axis background/.style={fill=white}
]
\addplot [color=mycolor1, style = solid , line width=1.0pt, mark=o, mark options={solid}, forget plot]
  table[row sep=crcr]{%
2	3	\\
3	3	\\
4	3	\\
5	3	\\
6	4	\\
7	5	\\
8	5	\\
9	5	\\
10	5	\\
11	5	\\
12	6	\\
13	6	\\
14	6	\\
};
\label{line:IC_d_dot}
\addplot [color=mycolor1, style = dashed, line width=1.0pt, mark=x, mark options={solid}, forget plot]
  table[row sep=crcr]{%
2	2	\\
3	2	\\
4	2	\\
5	2	\\
6	2	\\
7	2	\\
8	2	\\
9	2	\\
10	2	\\
11	2	\\
12	2	\\
13	2	\\
14	2	\\
};
\label{line:IC_d_s}
\addplot [color=mycolor2, style = solid, line width=1.0pt, mark=square, mark options={solid}, forget plot]
  table[row sep=crcr]{%
2	5	\\
3	6	\\
4	6	\\
5	6	\\
6	6	\\
7	6	\\
8	6	\\
9	6	\\
10	6	\\
11	6	\\
12	6	\\
13	6	\\
14	7	\\
};
\label{line:ICT_d_dot}
\addplot [color=black, style = solid, line width=1.0pt, mark=diamond, mark options={solid}, forget plot]
  table[row sep=crcr]{%
2	16	\\
3	16	\\
4	16	\\
5	16	\\
6	16	\\
7	16	\\
8	16	\\
9	16	\\
10	16	\\
11	16	\\
12	16	\\
13	16	\\
14	16	\\
};
\label{line:double}

\addplot [color=mycolor2, style = dashed, line width=1.0pt, mark=+, mark options={solid}, forget plot]
  table[row sep=crcr]{%
2	3	\\
3	3	\\
4	3	\\
5	3	\\
6	3	\\
7	3	\\
8	3	\\
9	3	\\
10	3	\\
11	3	\\
12	3	\\
13	3	\\
14	3	\\
};
\label{line:ICT_d_s}
\end{axis}

\begin{axis}[%
width=0.49\figurewidth,
height=\figureheight,
at={(0.5\figurewidth,0\figureheight)},
xlabel={level $J$},
ylabel = { number of IR iterations},
ylabel style={yshift=-1.2em},
xmin = 2,
xmax = 14,
ytick = {3,5,7,9,11},
yticklabels = {3,5,7,9,11},
xtick = {2,3,4,5,6,7,8,9,10,11,12,13,14},
xticklabels = {2,3,4,5,6,7,8,9,10,11,12,13,14},
yminorticks=true,
axis background/.style={fill=white}
]
\addplot [color=mycolor1, style = solid, line width=1.0pt, mark=triangle, mark options={solid}, forget plot]
  table[row sep=crcr]{%
2	11	\\
3	11	\\
4	11	\\
5	11	\\
6	10	\\
7	10	\\
8	10	\\
9	10	\\
10	10	\\
11	10	\\
12	9	\\
13	9	\\
14	9	\\
};
\label{line:IC_iter}
\addplot [color=mycolor2, style = solid, line width=1.0pt, mark=diamond, mark options={solid}, forget plot]
  table[row sep=crcr]{%
2	3	\\
3	3	\\
4	4	\\
5	4	\\
6	4	\\
7	4	\\
8	4	\\
9	4	\\
10	4	\\
11	4	\\
12	4	\\
13	4	\\
14	4	\\
};
\label{line:ICT_iter}

\end{axis}

\end{tikzpicture}
\caption{1D Poisson eq., FEM-P5 disc., solved by IR-V-cycle-IC. The plot on the left contains the values of $\dot{d}_{J,\mathrm{min}}$ and $d^{\solve}_{J,\mathrm{min}}$, i.e., the minimal values of $\dot{d}_J$ and $d^{\solve}_J$ such that the variant with $\dot{\varepsilon}_J=10^{-\dot{d}_J}$-precision and $\varepsilon^{\solve}_J=\varepsilon^{\store}_J=10^{-d^{\solve}_J}$-precision converges in the same number of IR iterations as the corresponding variant in double precision. The lines are labelled as
$\dot{d}_{J,\mathrm{min}}$ (\ref{line:IC_d_dot}), $d^{\solve}_{J,\mathrm{min}}$ (\ref{line:IC_d_s}) for the variant with IC(0)  and 
$\dot{d}_{J,\mathrm{min}}$ (\ref{line:ICT_d_dot}), $d^{\solve}_{J,\mathrm{min}}$ (\ref{line:ICT_d_s}) for the variant with 
ICT(dpt=$5\cdot10^{-3}$). For reference we also plot the number of digits for double precision (\ref{line:double}).
The plot on the right contains the number of IR iterations required for convergence for the variants in double precision with IC(0) (\ref{line:IC_iter}) and ICT(dpt=$5\cdot10^{-3}$) (\ref{line:ICT_iter}).}
\label{fig:1D_IC_results}
\end{figure}

\section{Numerical experiment: solving systems from high-order FE discretization of 3D elliptic PDEs}
\label{sec:num-exp-gingko}
The goal of this section is to show that the precisions used for applying the IC smoothing inside the V-cycle can be chosen to be lower than the precision used for computing the residual, prolongation, restriction, and correction, and that it may result in a significant reduction in the runtime and to energy and memory savings in comparison to the corresponding uniform double precision variant.

We solve systems of linear equations coming from the discretization of second order elliptic PDEs on a unit cube with homogeneous Dirichlet boundary conditions: find $u:(0,1)^3 \rightarrow \R$ such that 
\begin{equation*}
    -\nabla  (k\nabla u)   = 1 \ \textrm{in } (0,1)^3, \quad u = 0 \ \textrm{on } \partial (0,1)^3.
\end{equation*}
We consider two variants of the problem based on the coefficient function $k:(0,1)^3 \rightarrow \R$, a variant \texttt{Poisson} with $k \equiv 1$ and a variant \texttt{jump} with  
\begin{equation*}
    k(x) = 
    \begin{cases}
    1024, \quad x \in \left(0,\frac{1}{2} \right) \times \left(0,1 \right) \times \left( 0, 1 \right), \\
    1, \quad x \in \left(\frac{1}{2},1 \right) \times \left(0,1 \right) \times \left( 0, 1 \right).
    \end{cases}
\end{equation*}
The problems are discretized using the continuous Galerkin FE method with piecewise polynomials of degree five (FEM-P5) on a hierarchy of four uniformly refined meshes. The meshes are aligned with the plane segment where the jump in the coefficient function take place. The matrices are assembled in FEniCS as in the first experiment. We again use the scaling of the system and prolongation matrices described in \Cref{sec:scaling}. The numbers of DoFs on the individual levels are $6~859$, $59~319$, $493~039$ and $4~019~679$.

The iterative methods used to solve the algebraic systems on the finest level are implemented using the Ginkgo linear algebra library \cite{ginkgo,Ginkgo2024} and run on an NVIDIA GH200 GPU, with CUDA version 13.1.80, on the Kairos supercomputer at CALMIP, Toulouse. We chose the Ginkgo library for its support for implementing mixed precision algorithms.
The matrices, including the multigrid hierarchy, are generated in FEniCS in double precision and subsequently loaded on the GPU where the solution of the algebraic system using the Ginkgo library takes place. The matrices in Ginkgo are stored in the compressed sparse row (CSR) format with $32$-bit integer index type. The energy and memory measurements are obtained using the NVIDIA Management Library.
\subsection{IR-V-cycle-IC}
We solve the algebraic problems using the IR method with the inner solver chosen as the geometric V-cycle method (\Cref{alg:V-cycle}) with one iteration of IC(0) smoothing (\Cref{alg:ic_smoothing}) and with a direct solver based on Cholesky factorization on the coarsest level.
The methods are run with a zero initial approximate solution and stopped when the Euclidean norm of the relative residual is less than $10^{-10}$, i.e., $\| \vec{b} - \A \hat{\vec{x}}^{(i)} \| / \| \vec{b}\| \leq 10^{-10}$, where $\hat{\vec{x}}^{(i)}$ is the computed approximation.

The residual computation and the approximate solution update within IR are done in double precision.
We consider settings where the $\dot{\varepsilon}_j$-, $\varepsilon^{\store}_j$-, and $\varepsilon^{\solve}_j$- precisions do not vary on levels one to three, i.e., $\dot{\varepsilon}_j=\dot{\varepsilon}$, $\varepsilon^{\store}_j=\varepsilon^{\store}$, and $\varepsilon^{\solve}_j=\varepsilon^{\solve}$, $j=1,2,3$.
A summary of the mixed precision variants of the V-cycle-IC with concrete choices of the individual precisions can be found in \Cref{tab:3d-exp-precisions}.
``Single-half mix'' variant for the triangular solve means that the triangular solve algorithm uses single precision for the arithmetic operations and shared memory storage and half precision for global memory storage. A detailed description of the triangular solve algorithm can be found in \Cref{ape:triangular_solver}.
For the variants utilizing half precision, we use scaling of the right-hand-side vectors before applying the IC smoothing; see  \Cref{sec:scaling}. Vector scaling is also applied in the variant \texttt{h-s-h-sh-s} before calling the V-cycle and on levels $2$ and $1$ inside the V-cycle before projecting the residual vector to a coarser level.

\begin{table}
\caption{Summary of mixed precision variants of V-cycle-IC.}
\centering
\small
\label{tab:3d-exp-precisions}
\begin{tabular}{l|ccccc}
\toprule
variant & \makecell{residual, \\ prolongation,  \\ restriction, \\ correction\\ ($\dot{\varepsilon}$-prec.) } 
        & \makecell{IC \\ fact.} 
        & \makecell{$\matrx{L}_j$ stored in \\ ($\varepsilon^{\store}$-prec.)} 
        & \makecell{IC triang. \\ solve \\ ($\varepsilon^{\solve}$-prec.)}  & \makecell{coarsest-level \\ Cholesky \\ fact. and  solve } \\
\midrule
\texttt{d-d-d-d-d}  & double & double & double & double & double\\
\texttt{d-s-s-s-d} & double & single & single & single & double \\
\texttt{d-s-h-sh-d} & double & single & half & single-half mix & double \\
\texttt{s-s-s-s-s} & single & single & single & single & single \\
\texttt{s-s-h-sh-s} & single & single & half & single-half mix & single \\
\texttt{h-s-h-sh-s} & half & single & half & single-half mix & single \\
\bottomrule
\end{tabular}
\end{table}
\subsubsection*{Expectations based on the finite precision error analysis}

Since we are comparing variants of the IR-V-cycle with different choices of the precisions only inside the V-cycle, as in the previous experiment, it is enough to discuss how these precisions will affect the contraction factor of the V-cycle. We will use the results presented in \Cref{sec:IC-V-cycle}.
The approximate values of $\kappa^{\frac{1}{2}}_{A_j}$,  $\underline{\kappa}_{L_j}$, and $\|\matrx{L}_j^{-1} \|^2$ computed using the power method are summarized in \Cref{fig:cond_numbers_3d} for both the \texttt{Poisson} and \texttt{jump} problem.
We see that $\kappa^{\frac{1}{2}}_{A_j}$ grows by approximately a factor of two with each finer level while the values of $\underline{\kappa}_{L_j}$  and $\|\matrx{L}_j^{-1} \|^2$ stay nearly constant. The values of $\kappa^{\frac{1}{2}}_{A_j}$, $\underline{\kappa}_{L_j}$, and $\|\matrx{L}_j^{-1} \|^2$ corresponding to the \texttt{jump} problem are larger than the values corresponding to the \texttt{Poisson} problem.
We have also computed the following properties of the multigrid hierarchy, which are nearly the same for both problems and do not differ significantly on different levels:
$\| \A_j \| \approx 2.8$, $\||\A_j|\| \approx 3.8$, $m_{A_j} = 1331$, $\| \A^{-1}_{j-1}\|^{1/2} / \ \| \A^{-1}_j \|^{1/2} \approx 2 $, $\| \matrx{P}_j \| \approx 4.3$,  $\| | \matrx{P}_j | \| \approx 8$, $\underline{m}_{L_j} = 1331$. For $\underline{m}_{P_j}$ we have  $\underline{m}_{P_1}=3982$, $\underline{m}_{P_2}=8438$, $\underline{m}_{P_3}=8444$.
\begin{figure}
\centering
\setlength\figureheight{4.5cm}
\setlength\figurewidth{0.99\textwidth}
\definecolor{mycolor1}{RGB}{94, 201, 98}%
\definecolor{mycolor2}{RGB}{33, 145, 140}%
\definecolor{mred}{RGB}{237, 121, 83}%
\definecolor{mblue}{RGB}{65, 68, 135}%
\begin{tikzpicture}
\begin{axis}[%
width=0.4\figurewidth,
height=\figureheight,
at={(0\figurewidth,0\figureheight)},
title=\texttt{Poisson},
xlabel={level $j$},
ylabel = {},
ylabel style={yshift=-1.2em},
xmin = 1,
xmax = 3,
ymin = 1,
ymax = 100000,
ymode = log,
xtick = {0,1,2,3,4},
xticklabels = {0,1,2,3,4},
ytick = {1,10,100,1000,10000,100000},
yticklabels = {$1$,$10^1$,$10^2$,$10^3$,$10^4$,$10^5$},
yminorticks=true,
legend pos= outer north east,
legend style = {draw = none},
legend cell align=left,
axis background/.style={fill=white}
]
\addplot [color=mycolor1, style = solid, line width=1.0pt, mark=asterisk, mark options={solid}]
  table[row sep=crcr]{%
3	2.65E+02	\\
2	1.33E+02	\\
1	6.63E+01	\\
};
\addplot [color=mred, style = dashed, line width=1.0pt, mark=+, mark options={solid}]
  table[row sep=crcr]{%
3	9.09E+00	\\
2	9.02E+00	\\
1	8.85E+00	\\
};
\addplot [color=mblue, style = dashed, line width=1.0pt, mark=x, mark options={solid}]
  table[row sep=crcr]{%
3	2.84E+01	\\
2	2.84E+01	\\
1	2.81E+01	\\
};
\legend{$\kappa^{\frac{1}{2}}_{A_j}$, $\underline{\kappa}_{L_j}$, $\| \matrx{L}^{-1}_j \|^2$
}
\end{axis}

\begin{axis}[%
width=0.4\figurewidth,
height=\figureheight,
at={(0.5\figurewidth,0\figureheight)},
title=\texttt{jump},
xlabel={level $j$},
ylabel = {},
ylabel style={yshift=-1.2em},
xmin = 1,
xmax = 3,
ymin = 1,
ymax = 100000,
ymode = log,
xtick = {0,1,2,3,4},
xticklabels = {0,1,2,3,4},
ytick = {1,10,100,1000,10000,100000},
yticklabels = {$1$,$10^1$,$10^2$,$10^3$,$10^4$,$10^5$},
yminorticks=true,
legend pos= outer north east,
legend style = {draw = none},
legend cell align=left,
axis background/.style={fill=white}
]
\addplot [color=mycolor1, style = solid, line width=1.0pt, mark=asterisk, mark options={solid}]
  table[row sep=crcr]{%
1	1500.27027	\\
2	2999.801673	\\
3	5999.387222	\\
};
\addplot [color=mred, style = dashed, line width=1.0pt, mark=+, mark options={solid}]
  table[row sep=crcr]{%
3	2.91E+02	\\
2	2.88E+02	\\
1	2.81E+02	\\
};
\addplot [color=mblue, style = dashed, line width=1.0pt, mark=x, mark options={solid}]
  table[row sep=crcr]{%
3	2.92E+04	\\
2	2.89E+04	\\
1	2.82E+04	\\
};
\legend{$\kappa^{\frac{1}{2}}_{A_j}$, $\underline{\kappa}_{L_j}$,$\| \matrx{L}^{-1}_j \|^2$
}
\end{axis}
\end{tikzpicture}
\caption{3D \texttt{Poisson} and \texttt{jump} problems, FEM-P5 discretization. Properties of $\A_j$ and $\matrx{L}_j$ for IC(0).}
\label{fig:cond_numbers_3d}
\end{figure}

Using an analogous argument as for the previous experiment, with the difference that $\varepsilon^{\solve}$ and $\varepsilon^{\store}$ might differ, it can be shown that the $\A$-norm of the relative finite precision error after one V-cycle can be bounded as 
\begin{align*}
\frac{\|  \vec{y}_{\V} - \hat{\vec{y}}_{\V} \|_{\A}}{\| \vec{y} \|_{\A}} 
 & \lesssim 6   (\varepsilon^{\solve}  \underline{m}_{L_J,\varepsilon^{\solve}} + \varepsilon^{\store}) \| \A_J \| \| \matrx{L}^{-1}_J \|^2   \underline{\kappa}_{L_J} J 
 +  2 (C_{1,J} \| \matrx{L}^{-1}_J \|^2 + C_{2,J} ) \dot{\varepsilon} \kappa^{\frac{1}{2}}_{A_J}.
\end{align*}
We see that the requirements given by the estimate on the $\varepsilon^{\solve}$- and $\varepsilon^{\store}$-precision are lower than on the $\dot{\varepsilon}$-precision. The requirements on the $\varepsilon^{\store}$-precision are lower than the $\varepsilon^{\solve}$-precision.

\subsubsection*{Results}
We present time and energy measurements for the IC factorizations and Cholesky factorization in \Cref{tab:setup_ic_smoothing}, and time, energy, and memory measurements of the solve phase of the IR-V-cycle-IC method in \Cref{tab:res-IR-Vcycle-IC}. We do not measure the setup time of the multigrid hierarchy (generation of the system matrices and prolongation matrices), which is done in the FeniCS library in double precision.
The results are obtained by averaging measurements of a number of consecutive runs of the factorization or solve phases, respectively.

We first comment on the results for the IC and Cholesky factorizations measurements summarized in \Cref{tab:setup_ic_smoothing}.
We see that the time required for the IC and Cholesky factorizations is relatively small in comparison to the time of the solve phase of IR-V-cycle. The single precision variants are slightly faster and require less energy than the corresponding double precision variants. This is true for both the \texttt{Poisson} and \texttt{jump} problems.

We comment on the results of the IR-V-cycle solve phase measurements in \Cref{tab:res-IR-Vcycle-IC}.
For both of the problems, all variants besides the \texttt{h-s-h-sh-s} variant, which results in stagnation or an error, converge to the desired accuracy in the same number of IR iterations.
The fact that the solvers converge for both variants of the problem in the same number of iterations is probably due to the use of IC smoothing which is able to handle the difficulty caused by the jump in the diffusion coefficient $k$.
We see significant time, energy, and memory savings for the convergent mixed precision variants in comparison with the double precision variants. The improvements of the variants involving half precision are larger than that of the variants involving only single precision.
The variants \texttt{s-s-s-s-s} and \texttt{s-s-h-sh-s} are only slightly faster than the variants \texttt{d-s-s-s-d} and \texttt{d-s-h-sh-d}, respectively. This is likely because most of the time is spent in the triangular solves. Handling the other operations in single precision instead of in double precision thus has only a small effect.

The solve phase of the IR-V-cycle-IC solver is dominated by sparse matrix vector products and triangular solves, which are typically memory bound; see e.g., \cite{Grutzmacher2023}.
Converting values of an CSR matrix (with 32-bit integer index type) from double precision to single or half precision results in a matrix which requires approximately $0.67$ or $0.5$ of the original memory, respectively; see e.g., \cite[Section~6.2]{Tsai2024}.
The observed relative execution times of 0.86 (for variants with double and single precisions) and 0.72 (for variants with double, single, and half precisions) are within the expected range based on the reduction in memory.

\subsection{PCG-V-cycle-IC}
Since multigrid methods are frequently used as preconditioners, we also evaluate the performance of the preconditioned conjugate gradient method (PCG) with a geometric V-cycle preconditioner with one iteration of IC(0) for pre- and post-smoothing. The computations in PCG except for the application of the V-cycle-IC preconditioner are done in double precision. We consider the same mixed precision variants of the V-cycle-IC as for the IR method, i.e., the variant described in \Cref{tab:3d-exp-precisions}.
The methods are run with a zero initial approximate solution and stopped when the norm of the relative iteratively computed residual inside PCG is less than $10^{-10}$. We verified that in these experiments, 
the relative norm of the explicitly computed residual is also less than $10^{-10}$, i.e., $\| \vec{b} - \A \hat{\vec{x}}^{(i)} \| / \| \vec{b}\|\leq 10^{-10}$, where $\hat{\vec{x}}^{(i)}$ is the computed approximation.
The results are summarized in \Cref{tab:res-PCG-Vcycle-IC}. We see that the V-cycle with IC pre- and post- smoothing applied in low precision works well also as a preconditioner. The relative time, energy, and memory improvements of the mixed precision variants with respect to the double precision variant are analogous as for the IR method. All variants of the PCG method, however, converge in fewer iterations and require less time and energy than any of the variants of the IR method.

\begin{table}
\caption{3D \texttt{Poisson} and \texttt{jump} problems, FEM-P5 discretization. Time and energy measurement of IC factorization on levels $1$, $2$, and $3$, and of the Cholesky factorization on the coarsest-level in double and single precision. The relative time and energy are with respect to the double precision variant.}
\centering
\begin{tabular}{l|cccc|cccc}
\toprule
& \multicolumn{4}{c|}{\texttt{Poisson} } & \multicolumn{4}{c}{\texttt{jump} } \\
\midrule
precision & \makecell{time \\ $\left[\text{s} \right]$ }  & \makecell{rel. \\ time} & \makecell{ener. \\ $\left[ \text{kJ} \right]$ } & \makecell{rel. \\ ener.} & \makecell{time \\ $\left[\text{s} \right]$ }  & \makecell{rel. \\ time} & \makecell{ener. \\ $\left[ \text{kJ} \right]$ } & \makecell{rel. \\ ener.} \\
\midrule
\multicolumn{9}{c}{IC factorization sum of time and energy on levels 1, 2, and 3 } \\
\midrule
double  & 0.86 &  & 0.32 &  & 0.86 &  & 0.32 &  \\
single & 0.83 & 0.96 & 0.30 & 0.94 & 0.82 & 0.96 & 0.30 & 0.93 \\
\midrule
\multicolumn{9}{c}{Cholesky factorization on the coarsest level }\\
\midrule
double & 0.36 &  & 0.08 &  & 0.36 &  & 0.08 &  \\
single & 0.35 & 0.97 & 0.08 & 0.92 & 0.35 & 0.97 & 0.08 & 0.93 \\
\bottomrule
\end{tabular}
\label{tab:setup_ic_smoothing}
\end{table}
\begin{table}
\caption{3D \texttt{Poisson} and \texttt{jump} problems, FEM-P5 discretization. Results of the solve phase of IR-V-cycle-IC mixed precision variants. The relative time, energy, and memory are with respect to the double precision variant \texttt{d-d-d-d-d}. The best relative time, energy, and memory are highlighted in green and bold.}
\centering
\begin{tabular}{c|l|crccccc}
\toprule
& variant & \makecell{ IR \\ it.}   & \makecell{time \\ $\left[\text{s} \right]$ }  & \makecell{rel. \\ time} & \makecell{ener. \\ $\left[ \text{kJ} \right]$ } & \makecell{rel. \\ ener.}   & \makecell{mem. \\ $\left[ \text{GB} \right]$ } & \makecell{rel. \\ mem.} \\
\midrule
\parbox[t]{2mm}{\multirow{6}{*}{\rotatebox[origin=c]{90}{\texttt{Poisson}}}} 
& \texttt{d-d-d-d-d} &  49 & 14.38 &  & 5.83 &  & 17.5 &  \\
& \texttt{d-s-s-s-d} & 49 & 12.37 & 0.86 & 4.81 & 0.82 & 15.7 & 0.89 \\
& \texttt{d-s-h-sh-d} & 49 & 10.58 & 0.74 & 4.43 & 0.76 & 14.8 & \colorbox{_LimeGreen}{\textbf{0.84}} \\
& \texttt{s-s-s-s-s} & 49 & 12.30 & 0.86 & 4.71 & 0.81 & 16.7 & 0.95 \\
& \texttt{s-s-h-sh-s} & 49 & 10.51 & \colorbox{_LimeGreen}{\textbf{0.73}} & 4.38 & \colorbox{_LimeGreen}{\textbf{0.75}} & 15.9 & 0.90 \\
& \texttt{h-s-h-sh-s} & \multicolumn{5}{l}{stagnates at rel. res. $10^{-2}$} \\
\midrule
\parbox[t]{2mm}{\multirow{6}{*}{\rotatebox[origin=c]{90}{\texttt{jump}}}}& \texttt{d-d-d-d-d} & 49 & 14.39 &  & 5.90 &  & 17.5 &  \\
& \texttt{d-s-s-s-d} & 49 & 12.32 & 0.86 & 4.82 & 0.82 & 15.7 & 0.89 \\
& \texttt{d-s-h-sh-d} & 49 & 10.58 & 0.74 & 4.44 & 0.75 & 14.8 & \colorbox{_LimeGreen}{\textbf{0.84}} \\
& \texttt{s-s-s-s-s} & 49 & 12.27 & 0.85 & 4.71 & 0.80 & 16.7 & 0.95 \\
& \texttt{s-s-h-sh-s} & 49 & 10.52 & \colorbox{_LimeGreen}{\textbf{0.73}} & 4.37 & \colorbox{_LimeGreen}{\textbf{0.74}} & 15.8 & 0.90 \\
& \texttt{h-s-h-sh-s} &  \multicolumn{5}{l}{error} \\
\bottomrule
\end{tabular}
\label{tab:res-IR-Vcycle-IC}
\end{table}

\begin{table}
\caption{3D \texttt{Poisson} and \texttt{jump} problems, FEM-P5 discretization. Results of the solve phase of PCG-V-cycle-IC mixed precision variants. The relative time, energy and memory are with respect to the double precision variant \texttt{d-d-d-d-d}. The best relative time, energy and memory are highlighted in green and bold.}
\centering
\begin{tabular}{l|l|crccccc}
\toprule
& variant & \makecell{ PCG  \\ it.}   & \makecell{time \\ $\left[\text{s} \right]$ }  & \makecell{rel. \\ time} & \makecell{ener. \\ $\left[ \text{kJ} \right]$ } & \makecell{rel. \\ ener.}   & \makecell{mem. \\ $\left[ \text{GB} \right]$ } & \makecell{rel. \\ mem.} \\
\midrule
\parbox[t]{2mm}{\multirow{5}{*}{\rotatebox[origin=c]{90}{\texttt{Poisson}}}}
& \texttt{d-d-d-d-d} & 13 & 7.38 &  & 3.06 &  & 17.6 &  \\
& \texttt{d-s-s-s-d} & 13 & 6.29 & 0.85 & 2.50 & 0.82 & 15.7 & 0.89 \\
& \texttt{d-s-h-sh-d} & 13 & 5.35 & 0.73 & 2.26 & 0.74 & 14.8 & \colorbox{_LimeGreen}{\textbf{0.84}} \\
& \texttt{s-s-s-s-s} & 13 & 6.23 & 0.84 & 2.48 & 0.81 & 16.8 & 0.95 \\
& \texttt{s-s-h-sh-s} & 
13 & 5.31 & \colorbox{_LimeGreen}{\textbf{0.72}} & 2.25 & \colorbox{_LimeGreen}{\textbf{0.74}} & 15.9 & 0.90 \\
& \texttt{h-s-h-sh-s} &  \multicolumn{3}{l}{stagnates at rel. res. $10^{-1}$}      &  &  \\
\midrule
\parbox[t]{2mm}{\multirow{5}{*}{\rotatebox[origin=c]{90}{\texttt{jump}}}}
& \texttt{d-d-d-d-d} & 13 & 7.38 &  & 3.02 &  & 17.6 &  \\
& \texttt{d-s-s-s-d} & 13 & 6.28 & 0.85 & 2.45 & 0.81 & 15.8 & 0.90 \\
& \texttt{d-s-h-sh-d} & 13 & 5.36 & 0.73 & 2.25 & 0.75 & 14.8 & \colorbox{_LimeGreen}{\textbf{0.84}} \\
& \texttt{s-s-s-s-s} & 13 & 6.27 & 0.85 & 2.41 & 0.80 & 16.8 & 0.95 \\
& \texttt{s-s-h-sh-s} & 13 & 5.33 & \colorbox{_LimeGreen}{\textbf{0.72}} & 2.23 & \colorbox{_LimeGreen}{\textbf{0.74}} & 15.9 & 0.90 \\
& \texttt{h-s-h-sh-s} & \multicolumn{3}{l}{error}      \\
\bottomrule
\end{tabular}
\label{tab:res-PCG-Vcycle-IC}
\end{table}

\section{Conclusion}

We present a mixed precision formulation of the multigrid \mbox{V-cycle} method with general assumptions on the finite precision errors of the coarsest-level solver and smoothers. We derive a bound on the relative finite precision error of the V-cycle which gives insight into how the finite precision errors from the individual components of the V-cycle may affect the overall finite precision error.
The presented approach enables analyses of V-cycles with various (mixed precision) coarsest-level solvers and smoothers. 

In this work, we focus on mixed precision smoothers based on IC factorization. 
We present theoretical results implying that in certain settings the precisions used for applying the IC smoothing could be significantly lower than the precision used for computing the residual, restriction, prolongation, and correction on a given level. Experiments on GPUs using the Ginkgo library show a significant speedup and energy savings of variants with low precisions IC smoothing.

In future work, we would like to apply these results to analyze mixed precision multigrid methods with other frequently used smoothers, e.g., Gauss-Seidel or SOR smoothers.
Moreover, it would be interesting to investigate if the presented results can be extended to methods with both pre- and post- smoothing and subsequently used to study the effects of finite precision errors in PCG with multigrid preconditioning.

\bibliographystyle{siamplain}
\bibliography{references}
\appendix
\section{Derivation of inequalities  \eqref{eq:A_norm_by_Euclid} - \eqref{eq:A^-1_v_by_Euclid}}
\label{ape:norms}

\noindent Ad. \eqref{eq:A_norm_by_Euclid}: 
 $\| \vec{v} \|^2_\A  = \langle \A \vec{v} , \vec{v} \rangle \leq \| \A \| \langle \vec{v} , \vec{v} \rangle = \| \A \| \| \vec{v} \|^2. $
\begin{align*}
&\text{Ad. \eqref{eq:Euclid_by_A_norm}: }  
     \| \vec{v} \|^2  = \langle \vec{v},\vec{v}  \rangle 
    = \langle \A^{-\frac{1}{2}}\A^{\frac{1}{2}} \vec{v}, \A^{-\frac{1}{2}}\A^{\frac{1}{2}} \vec{v}  \rangle
    = \langle \A^{-\frac{1}{2}}\A^{-\frac{1}{2}}\A^{\frac{1}{2}} \vec{v}, \A^{\frac{1}{2}} \vec{v}  \rangle \\
    & \qquad \qquad \qquad \quad = \langle \A^{-1} \A^{\frac{1}{2}} \vec{v}, \A^{\frac{1}{2}} \vec{v}  \rangle 
    \leq \| \A^{-1} \| \langle \A^{\frac{1}{2}}\vec{v},\A^{\frac{1}{2}}\vec{v}  \rangle \\ 
    & \qquad \qquad \qquad \quad =\| \A^{-1}  \| \langle \A\vec{v},\vec{v}  \rangle 
    = \| \A^{-1} \| \| \vec{v}  \|^2_{\A}. \\ 
&\text{Ad. \eqref{eq:A*v_by_A_norm}: } 
   \| \A\vec{v} \|^2  = \langle \A \vec{v}, \A\vec{v}  \rangle
    = \langle \A^{\frac{1}{2}} \A^{\frac{1}{2}} \vec{v}, \A^{\frac{1}{2}} \A^{\frac{1}{2}} \vec{v}  \rangle =  \langle \A^{\frac{1}{2}}\A^{\frac{1}{2}} \A^{\frac{1}{2}} \vec{v},  \A^{\frac{1}{2}} \vec{v}  \rangle \\ 
    & \qquad \qquad \qquad \qquad =  \langle \A \A^{\frac{1}{2}} \vec{v},  \A^{\frac{1}{2}} \vec{v}  \rangle \leq \| \A \| \langle \A^{\frac{1}{2}} \vec{v},  \A^{\frac{1}{2}} \vec{v}  \rangle  = \| \A \| \langle \A \vec{v} ,\vec{v}  \rangle = \| \A \| \| \vec{v}  \|^2_{\A}.\\   
&\text{Ad. \eqref{eq:A^-1_v_by_Euclid}: } 
     \| \A^{-1} \vec{v} \|^2_{\A}   = \langle \A \A^{-1} \vec{v},\A^{-1}\vec{v} \rangle
    = \langle  \vec{v},\A^{-1}\vec{v}
    \rangle     
    \leq \| \A^{-1} \| 
 \langle  \vec{v} ,\vec{v}  \rangle = \| \A^{-1} \| \| \vec{v}  \|^2.
\end{align*}
\section{Derivation of multigrid related inequalities \eqref{eq:special_projection_A_Cnorm_bound} - \eqref{eq:v[2]_C_A_C_norm_bound}}
\label{ape:mg}
Ad. \eqref{eq:special_projection_A_Cnorm_bound}: We rewrite $\| \A^{-1}_{\C} \matrx{P}^{\top} \vec{v} \|^2_{\A_{\C}}  $ as 
\begin{align}
\nonumber
\| \A^{-1}_{\C} \matrx{P}^{\top} \vec{v} \|^2_{\A_{\C}}  
&= \langle  \A_{\C}\A^{-1}_{\C} \matrx{P}^{\top} \vec{v} , \A^{-1}_{\C} 
\nonumber
\matrx{P}^{\top} \vec{v}  \rangle  
=\langle  \matrx{P}^{\top} \vec{v} , \A^{-1}_{\C} \matrx{P}^{\top} \vec{v}  \rangle
 = \langle  \vec{v} , \matrx{P} \A^{-1}_{\C} \matrx{P}^{\top} \vec{v}  \rangle \\
&=\langle  \A^{\frac{1}{2}} \A^{-\frac{1}{2}}\vec{v} , \matrx{P} \A^{-1}_{\C} \matrx{P}^{\top} \A^{\frac{1}{2}}\A^{-\frac{1}{2}}\vec{v}  \rangle 
\label{eq:special_projection_A_Cnorm_bound_rewrite}
=\langle  \A^{-\frac{1}{2}}\vec{v} , \A^{\frac{1}{2}} \matrx{P} \A^{-1}_{\C} \matrx{P}^{\top} \A^{\frac{1}{2}}\A^{-\frac{1}{2}}\vec{v}  \rangle. 
\end{align}
Since
$\A^{\frac{1}{2}} \matrx{P} \A^{-1}_{\C} \matrx{P}^{\top} \A^{\frac{1}{2}}  = \A^{\frac{1}{2}} \matrx{P} (\matrx{P}^{\top} \A \matrx{P})^{-1} \matrx{P}^{\top} \A^{\frac{1}{2}} =\A^{\frac{1}{2}} \matrx{P} ( ( \A^{\frac{1}{2}} \matrx{P} )^{\top} \A^{\frac{1}{2}} \matrx{P} )^{-1}  ( \A^{\frac{1}{2}} \matrx{P} )^{\top}   $ 
is the orthogonal projection onto the range of $\A^{\frac{1}{2}} \matrx{P}$,  
$\|  \A^{\frac{1}{2}} \matrx{P} \A^{-1}_{\C} \matrx{P}^{\top} \A^{\frac{1}{2}} \|\leq 1$.     
Combining this and \eqref{eq:special_projection_A_Cnorm_bound_rewrite} leads to
\begin{equation*}
\| \A^{-1}_{\C} \matrx{P}^{\top} \vec{v} \|^2_{\A_{\C}}
\leq \langle  \A^{-\frac{1}{2}}\vec{v} , \A^{-\frac{1}{2}}\vec{v}  \rangle 
= \langle   \vec{v} ,\A^{-1} \vec{v}  \rangle 
= \langle   \A \A^{-1} \vec{v} ,\A^{-1} \vec{v}  \rangle 
 = \| \A^{-1} \vec{v} \|^2_{\A}.
\end{equation*}

\noindent Ad. \eqref{eq:B_C_Anorm_bound}: Using assum. \eqref{eq:perturbed_coarse_solve_A_norm},
$\| \matrx{M}_{\C}  \A_{\C} \|_{\A_{\C}} \leq \| \matrx{I}_{\C} \|_{\A_{\C}} + \|\matrx{I}_{\C} -  \matrx{M}_{\C} \A_{\C}  \|_{\A_{\C}} < 2 .$

\noindent Ad. \eqref{eq:v[4]_Anorm_bound}: Using assum. \eqref{eq:TG_exact_rho},
$\| \vec{v}^{[4]} \|_{\A} = \| \vec{y}_{\TG} \|_{\A} \leq \| \vec{y}_{\TG}  - \vec{y}  \|_{\A} + \| \vec{y}\|_\A\leq 2 \| \vec{y}\|_\A.$

\noindent Ad. \eqref{eq:A-1_r^[1]_bound}: Using $\A \vec{y} = \vec{f}$, $\vec{r}^{[1]}=\vec{f} - \A \matrx{M} \vec{f}$, and the assumption \eqref{eq:smoother_assumption_Anorm} results in
\begin{equation*}
\| \A^{-1}  \vec{r}^{[1]}\|_{\A} 
=  \| \A^{-1}  (\vec{f} - \A \matrx{M} \vec{f})\|_{\A}
 =  \|  \vec{y} -  \matrx{M} \A \vec{y}\|_{\A}
 \leq  \| \matrx{I} - \matrx{M} \A\|_{\A} \| \vec{y} \|_{\A} 
\leq    \| \vec{y} \|_{\A}.    
\end{equation*}

\noindent Ad. \eqref{eq:r^[1]_bound}: Using \eqref{eq:A-1_r^[1]_bound} yields
$\|  \vec{r}^{[1]}\| \leq \| \A\|^{\frac{1}{2}} \| \A^{-1}  \vec{r}^{[1]}\|_{\A} \leq \| \A\|^{\frac{1}{2}} \| \vec{y} \|_{\A}$. 

\noindent Ad. \eqref{eq:A-1_Cr[2]_C_bound}: Using $\vec{r}^{[2]}_{\C} = \matrx{P}^{\top} \vec{r}^{[1]}$, \eqref{eq:special_projection_A_Cnorm_bound} and \eqref{eq:A-1_r^[1]_bound} results in
\begin{equation*}
\| \A^{-1}_\C \vec{r}^{[2]}_\C \|_{\A_\C} = \| \A^{-1}_\C \matrx{P}^{\top} \vec{r}^{[1]}\|_{\A_\C}  \leq 
\| \A^{-1} \vec{r}^{[1]} \|_{\A} \leq 
\| \vec{y}\|_\A.
\end{equation*}

\noindent Ad. \eqref{eq:v[2]_C_A_C_norm_bound}: Using  $\vec{v}^{[2]}_{\C}=\matrx{M}_{\C} \vec{r}^{[2]}_{\C}$ and bounds \eqref{eq:B_C_Anorm_bound} and \eqref{eq:A-1_Cr[2]_C_bound} results in
\begin{equation*}
\| \vec{v}^{[2]}_{\C} \|_{\A_{\C}} 
 = \| \matrx{M}_{\C} \A_{\C}  \A_{\C}^{-1} \vec{r}^{[2]}_{\C} \|_{\A_{\C}} 
\leq \| \matrx{M}_{\C} \A_{\C} \|_{\A_{\C}}  \| \A_{\C}^{-1} \vec{r}^{[2]}_{\C} \|_{\A_{\C}} 
\leq  2  \| \vec{y} \|_{\A}. 
\end{equation*}

\section{Proof of \Cref{lemma:substitution}}
\label{ape:sub_fp}
We present a proof for a lower-triangular matrix $\matrx{T}$. The proof for an upper-triangular matrix $\matrx{T}$ is analogous. 
We will use the following lemma.
\begin{lemma}{\cite[Lemma~8.4]{Higham2002}}\label{lemma:Higham-L8.4}
Let $k$ be a natural number an let $\delta$, $\alpha_i$, $i=1,\ldots,k-1$, $\beta_i$, $i=1,\ldots,k$, be numbers belonging to a finite precision arithmetic with unit roundoff $\varepsilon$ and $k\varepsilon<1$. Computing
  $  \gamma = ( \delta -\sum^{k-1}_{i=1}\alpha_i \beta_i)/\beta_k$
in $\varepsilon$-precision results in $\hat{\gamma}$ satisfying, no matter the order of evaluation of the sum,
\begin{equation*}
\beta_k(1+\theta^{(0)}_k)\hat{\gamma} = \delta-\sum^{k-1}_{i=1}\alpha_i \beta_i (1+\theta^{(i)}_k),
\quad \text{where} \quad |\theta^{(i)}_k|\leq \frac{k \varepsilon} {1-k\varepsilon},  \ i=0,1,\ldots,k.
\end{equation*}
\end{lemma}

We use induction on the size of the leading sub-matrices. Let $t_{i,j}$ and $e_{i,j}$ denote the entries of matrices $\matrx{T}$ and $\matrx{E}$, respectively, in the $i$th row and $j$th column, and let $b_i$ denote the $i$th entry of the vector $\vec{b}$.

We start by showing that the statement holds for the leading sub-matrix of size $1\times1$.
Using \Cref{lemma:Higham-L8.4} for $k=1$, computing $x_1 = b_1/t_{1,1}$ in $\varepsilon$-precision results in $\hat{x}_1$ satisfying $t_{1,1}(1+\theta^{(0)}_1)\hat{x}_1 = b_1$, where $|\theta^{(0)}_1|\leq \varepsilon/(1-\varepsilon)$. We can take $e_{1,1} = t_{1,1}\theta^{(0)}_1$.

Assume that the statement holds for the leading sub-matrix of size $ n \times n$.
We will show that it holds also for the leading sub-matrix of size $(n+1)\times(n+1)$. 
Using the induction assumption and the fact that $\matrx{T}$ is a lower triangular matrix, it only remains to show the existence of suitable entries in the $(n+1)$th row of $\matrx{E}$.
Let $\hat{x}_i$, $i=1,\ldots,n$ denote the computed entries of $\hat{\vec{x}}$ after $n$ steps of the substitution.
The ($n+1)$th substitution step consists of computing
\begin{equation*}
    x_{n+1} = (b_{n+1} - \sum^{n}_{i=1} \hat{x}_i t_{n+1,i} )/t_{n+1,n+1}.
\end{equation*}
Since we assume there is a maximum of $m_T$ nonzero elements in a row of $\matrx{T}$, the sum consists of a maximum of $m_T-1$ nonzero terms. The equation can be rewritten as
\begin{equation*}
    x_{n+1} = (b_{n+1} - \sum_{\ell; t_{n+1,\ell}\neq 0 } \hat{x}_\ell t_{n+1,\ell} )/t_{n+1,n+1}.
\end{equation*}
Using \Cref{lemma:Higham-L8.4} in this setting yields
\begin{equation*}
t_{n+1,n+1}(1+\theta^{(0)}_{m_T}) \hat{x}_{n+1} = b_{n+1} -  \sum_{\ell; t_{n+1,\ell}\neq 0 } \hat{x}_\ell t_{n+1,\ell}(1+\theta^{(\ell)}_{m_T}),
\end{equation*}
where $|\theta^{(i)}_{m_T}|\leq (m_T \varepsilon)/(1-m_T \varepsilon)$, $i=0,\ldots,m_T-1$. Taking $e_{n+1,\ell} =  t_{n+1,\ell} \theta^{(\ell)}_{m_T} $, for $\ell$ such that $t_{n+1,\ell}\neq 0$ and $e_{n+1,n+1} =  t_{n+1,n+1}\theta^{(0)}_{m_T}$ implies that the statement holds also for the $(n+1)\times(n+1)$ leading sub-matrix.

\section{Triangular solve algorithm}
\label{ape:triangular_solver}
We describe the details of the mixed precision triangular solve algorithm used in the experiments in \Cref{sec:num-exp-gingko}.
We use the CUDA \emph{syncfree} variant of the Ginkgo triangular solve. Its  sketch can be found in \Cref{alg:triangular_solve}.

In this algorithm, one thread is assigned to one row of the matrix. The threads are organized in groups sharing access to fast shared memory.
Each thread waits until a result associated with an above row that is necessary for its computation is available. This result is read either from the shared or the global memory depending on whether the threads are part of the same group.

We consider the case where the $\varepsilon_{\text{shared}}$-precision in which values are stored in the shared memory and the $\varepsilon_{\text{global}}$-precision used for storing values in the global memory may differ. All floating point operations in the algorithm are done in the  $\varepsilon_{\text{shared}}$-precision.

We apply this algorithm either with uniform precision configurations, where $\varepsilon_{\text{shared}}$- and $\varepsilon_{\text{global}}$- precisions are both chosen as double or single or in the single-half mix variant, where $\varepsilon_{\text{shared}}$-precision is chosen as single and $\varepsilon_{\text{global}}$-precision as half.

\begin{algorithm}
\caption{Lower triangular solve algorithm}
\label{alg:triangular_solve}
\begin{algorithmic}[]
\STATE{Inputs: lower triangular matrix $\matrx{T}$ with values stored in $\varepsilon_{\text{global}}$-precision and the right-hand side vector $\vec{b}$ in $\varepsilon_{\text{global}}$-precision.}
\STATE{Output: the solution vector $\vec{x}$ in $\varepsilon_{\text{global}}$-precision.} 
\STATE{}
\STATE{$i $ \COMMENT{row index = thread index; one thread is assigned to one row}} 
\STATE{$s = 0$}
\FOR{ $j=1$; $(j < i)\text{ and }(t_{i,j} \neq 0)$; $j++$  }
\IF{the thread assigned to the $j$th row is in the same group}
\STATE{wait until the $j$th row result $r_j$ is available in shared memory}
\STATE{$ s = s + t_{i,j} \cdot r_j $ \COMMENT{Compute the update of $s$ in $\varepsilon_{\text{shared}}$-precision}}
\ELSE
\STATE{wait until the $j$th row result $x_{j}$ is available in global memory}
\STATE{$ s = s + t_{i,j} \cdot x_{j}$ \COMMENT{Compute the update of $s$ in $\varepsilon_{\text{shared}}$-precision}}
\ENDIF
\ENDFOR
\STATE{$r = (b_{i} - s) / t_{i,i}$ \COMMENT{Compute the row result in $\varepsilon_{\text{shared}}$-precision}} 
\STATE{store $r$ in shared memory in $r_{j}$ in $\varepsilon_{\text{shared}}$-precision} 
\STATE{store $r$ in  global memory in $x_{j}$ in $\varepsilon_{\text{global}}$-precision}
\end{algorithmic}
\end{algorithm}

\end{document}